\documentclass[12pt]{article}
\title{Birational geometry and deformations of 
nilpotent orbits}
\author{Yoshinori Namikawa}
\date{ }
\usepackage{amsmath,amssymb,amsthm, amscd}
\chardef\bslash=`\\
\newtheorem{Thm}{Theorem}[section]
\newtheorem{Cor}[Thm]{Corollary}
\newtheorem{Lem}[Thm]{Lemma}
\newtheorem{Prop}[Thm]{Proposition}
\newtheorem{Def}{Definition}
\newtheorem{Rque}[Thm]{Remark}

\newtheorem{Exam}[Thm]{Example}

\def\0{{\mathcal O}}
\def\g{{\mathfrak g}}
\def\q{{\mathfrak q}}
\def\p{{\mathfrak p}}
\def\h{{\mathfrak h}}
\def\k{{\mathfrak k}}
\def\l{{\mathfrak l}}
\def\b{{\mathfrak b}}

\def\z{{\mathfrak z}}

\begin{document}
\maketitle

Let $G$ and $\g$ be a complex simple Lie group and 
its Lie algebra. An orbit $\mathcal{O} \subset \g$ for 
the adjoint action of $G$, is called a nilpotent orbit if 
it consists of nilpotent elements. 
For a parabolic subgroup $P$ of $G$, there exists a 
unique nilpotent orbit $\mathcal{O}$ such that 
$\mathcal{O} \cap n(\p)$ is an open dense subset of 
$n(\p)$,  
where $n(\p)$ is the nilradical of $\p$. Then this 
orbit $\mathcal{O}$ is called the Richardson orbit for $P$. 
The closure $\bar{\mathcal{O}}$ of $\mathcal{O}$ 
admits symplectic singularities and it is an interesting 
object in the birational geometry and the singularity 
theory. There is a generically finite, projective morphism 
(Springer map) from the cotangent bundle of $G/P$ to 
$\bar{\mathcal{O}}$:  
$$ s: T^*(G/P) \to \bar{\mathcal{O}}. $$ 
Let $T^*(G/P) \stackrel{\pi}\to \tilde{\mathcal{O}} 
\to \bar{\mathcal{O}}$ be the Stein factorization 
of $s$. 
When $\mathrm{deg}(s) = 1$, $s$ is a crepant 
resolution of $\bar{\mathcal{O}}$. But, in this case, 
$\bar{\mathcal{O}}$ generally has other crepant resolutions. 
In [Na], we have shown that all such crepant resolutions 
are obtained as the Springer resolutions $s'$ for certain 
(finitely many)  parabolic 
subgroups $P'$, and classified the conjugacy classes of 
such parabolic subgroups in terms of  Dynkin diagrams.  
The (closed) movable cone $\overline{\mathrm{Mov}}(s) 
\subset \mathrm{Pic}(T^*(G/P))\otimes {\mathbf{R}}$ (cf. (P.2) below) 
is decomposed by the nef cones $\overline{\mathrm{Amp}}(s')$ 
of various Springer resolutions $s'$.  
In general, $\overline{\mathrm{Mov}}(s) \ne \mathrm{Pic}(T^*(G/P))
\otimes \mathbf{R}$. The first purpose of this paper, is to 
understand all parts of $\mathrm{Pic}(T^*(G/P))\otimes \mathbf{R}$ 
in terms of ample cones of certain varieties. 
The second purpose is, to prove a similar result to [Na] when 
the Springer map $s$ has $\mathrm{deg}(s) > 1$.   
  
Our strategy is to use  Brieskorn-Slodowy diagrams. 
Let us fix a maximal torus $T$ of $P$, 
and denote by $\h$ its Lie algebra. Let $l(\p)$ be the 
Levi factor of $\p$ containing $\h$, and $\k(\p)$ the centralizer 
of $l(\p)$.  Let $r(\p)$ (resp. $n(\p)$) be the solvable radical 
(resp. nilradical) of $\p$.  The Brieskorn-Slodowy diagram 
is the following commutative diagram (for details, see (P.1)): 
\newpage 
  
$$ G \times^P r(\p) \to Gr(\p) $$ 
$$ \downarrow \hspace{1.0cm} \downarrow $$ 
$$ \mathfrak{k}(\p) \to \mathfrak{h}/W. $$
The variety $X_{\p} := G \times^P r(\p)$ is a flat 
deformation of the cotangent bundle $T^*(G/P) = 
G \times^P n(\p)$ over the parameter space $\k(\p)$ 
with some universal property; this is indeed the universal 
Poisson deformation of the complex symplectic variety 
$T^*(G/P)$.  In turn, the fiber product 
$\mathfrak{k}(\p) \times_{\mathfrak{h}/W} 
Gr(\p)$ has the projection map to 
$\k(\p)$. This map is equi-dimensional. The 
central fiber is the closure $\bar{\mathcal{O}}$ 
of the Richardson orbit 
${\mathcal{O}}$, and the fiber over a general point 
$t \in \k(\p)$ is the semi-simple orbit containing $t$. 
Lemma 1.1 claims that the induced map 
$$ \mu'_{\p}: X_{\p} \to \mathfrak{k}(\p) \times_{\mathfrak{h}/W} 
Gr(\p)$$ is a birational projective 
morphism, and it is an isomorphism over a general 
point $t \in \k(\p)$. 
Let  $$X_{\p} \stackrel{\mu_{\p}}\to Y_{\k(\p)} \to     
\mathfrak{k}(\p) \times_{\mathfrak{h}/W} 
Gr(\p)$$ be the Stein factorization of $\mu'_{\p}$. 
In other words, $Y_{\k(\p)}$ is the normalization of 
$\mathfrak{k}(\p) \times_{\mathfrak{h}/W} 
Gr(\p)$. 
Here we must remark that $Gr(\p)$ actually 
depends only on $\k(\p)$ (or $l(\p)$) (cf. Theorem 0.1). 
So, $Y_{\k(\p)}$ depends only on $\k(\p)$ as this notation indicates.   
As we remark in Observation 1, the central fiber $Y_{\k(\p),0}$ of 
$Y_{\k(\p)} \to \k(\p)$, has $\tilde{\mathcal{O}}$ as its 
normalization. This is an important point; when $\mathrm{deg}(s) 
> 1$, the natural map $Y_{\k(\p),0} \to \bar{\mathcal{O}}$ 
has degree $> 1$.  Since $Y_{\k(\p),0}$ coincides with $\tilde{\mathcal{O}}$ 
at least set theoretically (cf. Observation 1), 
$Y_{\k(\p)} \to \k(\p)$ is {\em almost} a deformation of 
$\tilde{\mathcal{O}}$. By Lemma 1.2, $\mu_{\p}$ is a 
crepant resolution which is an isomorphism in codimension one. 
Now let us fix one parabolic subalgebra $\p_0$ containing 
$\h$, and put $l_0 := l(\p_0)$ and $\k_0 := \k(\p_0)$. 
Define $\mathcal{S}(l_0)$ as the set of all 
parabolic subalgebras $\p$ containing 
$\h$ such that {\em their Levi factors $l(\p)$ coincide with 
$l_0$}.   
\vspace{0.2cm}
 
Our main theorem (= Theorem 1.3) states that 
any crepant resolution of $Y_{\k_0}$ has the form 
$\mu_{\p}: X_{\p} \to Y_{\k_0}$ with 
$\p \in \mathcal{S}(l_0)$. Note that $\k(\p) = \k_0$ 
for $\p \in \mathcal{S}(l_0)$.       
Moreover, $\mu_{\p} \ne \mu_{\p'}$ if $\p \ne \p'$. 
\vspace{0.2cm}

A main ingredient in the proof of Theorem 1.3 is the 
notion of a {\em twist}.  Each parabolic subalgebra $\p$ 
corresponds to a marked Dynkin diagram when a Borel 
subalgebra $\b \subset \p$ is fixed (cf. (P.1)). 
Let $D$ be such a marked Dynkin diagram determined by 
$\p \in \mathcal{S}(l_0)$. Some vertices of $D$ are 
marked with black. For each marked vertex $v$, we can make  
a new parabolic subalgebra $\p' \in \mathcal{S}(l_0)$.  
The new parabolic subalgebra $\p'$ is called the twist of $\p$ by  
$v$.  Depending on $v$, the twists are classified into two 
classes (twists of the 1-st kind, and twists of the 2-nd kind). 
If the twist is of the 1-st kind, then $\p'$ is not conjugate 
to $\p$. But, for a twist of the 2-nd kind, $\p'$ is 
conjugate to $\p$. The two varieties $X_{\p}$ and $X_{\p'}$ 
are related by a flop.  The proof 
of Theorem 1.3 goes as follows. Let $\mu: X \to Y_{\k_0}$  
be an arbitrary crepant resolution. Take a $\mu$-ample line 
bundle $L$ on $X$ and let $L^{(0)} \in \mathrm{Pic}(X_{\p_0})$ 
be its proper transform. If $L^{(0)}$ is not $\mu_{\p_0}$-nef, we have  
an extremal birational contraction map $g_0: X_{\p_0} \to \bar{X}_{\p_0}$ 
with respect to $L^{(0)}$. By the definition of $g_0$, $L^{(0)}$ is  
$g_0$-negative. Moreover, since $\mu_{\p_0}$ is an isomorphism in codimension 
one (cf. Lemma 1.2), $g_0$ is an isomorphism in codimension one. A flop of $g_0$ is a 
diagram $X_{\p_0} \stackrel{g_0}\to \bar{X}_{\p_0} \stackrel{g'_0}\leftarrow X_1$ 
such that the proper transform $L^{(1)} \in \mathrm{Pic}(X_1)$ of $L^{(0)}$ 
is $g'_0$-ample (cf. (P.2)).  In our case, $X_1$ coincides with $X_{\p_1}$ 
for a twist $\p_1$ of $\p_0$.  If  $L^{(1)}$ is not $\mu_{\p_1}$-nef, then 
we repeat this process. In this way,  we continue the flops $X_{\p_0} --\to X_{\p_1} --\to ...$ 
as long as  the proper transform of $L^{(0)}$ is not nef over $Y_{\k_0}$.  But, since  
$\mathcal{S}(l_0)$ is a finite set, the sequence terminates at some 
$X_{\p_k}$. Then the proper transform $L^{(k)} \in \mathrm{Pic}(X_{\p_k})$ 
of $L^{(0)}$ is $\mu_{\p_k}$-nef, and $X = X_{\p_k}$(cf. Corollary 1.5).    

For a parabolic subgroup $P$ with $\p \in \mathcal{S}(l_0)$, 
$\mathrm{Pic}(G/P)_{\mathbf{R}} := 
\mathrm{Pic}(G/P) \otimes{\mathbf R}$ is canonically 
identified with $M(L_0)_{\mathbf{R}}$ (cf. (P.3)), where 
$$M(L_0) :=   
\mathrm{Hom}_{\mathrm{alg}.\mathrm{gp}.}(L_0, \mathbf{C}^*).$$ 
On the other hand, $\mathrm{Pic}(G/P)$ and $\mathrm{Pic}(X_{\p})$ 
are identified by the projection map $X_{\p} \to G/P$. 
We therefore have an isomorphism 
$$ \Phi_P: \mathrm{Pic}(X_{\p})_{\mathbf{R}} \to M(L_0)_{\mathbf{R}}. $$
Let us take another $P'$ with a fixed Levi part $L_0$.  
The isomorphism above is natural in the sense that 
$(\Phi_{P'})^{-1}\circ \Phi_P : \mathrm{Pic}(X_{\p})_{\mathbf{R}} 
\cong \mathrm{Pic}(X_{\p'})_{\mathbf{R}}$ coincides with the  
isomorphism defined by the proper transform by  
$X_{\p} --\to X_{\p'}$ (Observation 2).  
Moreover, the image $\Phi_P(\overline{\mathrm{Amp}}
(\mu_{\p}))$ of the nef cone for $\mu_{\p}$ by $\Phi_P$, 
can be described explicitly in terms of dominant characters.   
An important corollary to (the proof of) Theorem 1.3 is 
Remark 1.6, which says that  
$$M(L_0)_{\mathbf{R}} = 
\bigcup_{\p \in \mathcal{S}(l_0)}\Phi_P(\overline{\mathrm{Amp}}
(\mu_{\p})). $$ On the $M(L_0)_{\mathbf{R}}$, a subgroup of 
the Weyl group $W$ of $\g$ acts. More explicitly, this 
subgroup is the normalizer $N_W(L_0)$ of $L_0$. 
Let $W(L_0)$ be the Weyl group of $L_0$. Then 
$N_W(L_0)/W(L_0)$ acts effectively on $M(L_0)_{\mathbf{R}}$. 
In Section 2, we shall prove that the set $\mathcal{S}(l_0)$ 
contains exactly $N \cdot \sharp(N_W(L_0)/W(L_0))$ elements, where 
$N$ is the number of the conjugacy classes of parabolic subalgebras 
contained in $\mathcal{S}(l_0)$. 
The number $N$ can be calculated explicitly in terms of marked Dynkin 
diagrams [Na], Definition 1 (see also [Ri]).   
Howlett [Ho] gives an explicit description of $N_W(L_0)/W(L_0)$.   
When the Springer map $s: T^*(G/P_0) \to \bar{\mathcal{O}}$ 
is a resolution\footnote{This condition can be translated to 
a combinatorial condition for $L_0$ when $G$ is classical [He].}, 
the movable cone $\mathrm{Mov}(\mu_{\p})$ 
for $\mu_{\p}$ is a fundamental domain of the $N_W(L_0)/W(L_0)$-action (Proposition 2.3).   
In Section 3, we shall study the crepant resolutions of $\tilde{\mathcal O}$.   
When $\mathrm{deg}(s) = 1$, this is [Na], Theorem 6.1.  The problem is actually 
in the case where $\mathrm{deg}(s) > 1$.  Now we can use the map $\mu_{\p}: 
X_{\p} \to Y_{\k(\p)}$ to study the map 
$\pi: T^*(G/P) \to \tilde{\mathcal{O}}$, because $\tilde{\mathcal O}$ is the normalization 
of $Y_{\k(\p),0}$ and $\pi$ is the Stein factorization of the birational map  
$\mu_{\p,0}$.  
The main result in this section is Corollary 3.4. 
In order to state the result, we need to go back to the notion ``twist" again. 
As remarked above, there are two kind of twists. In Section 3, we divide 
the twists of the 2-nd kind into two classes: small twists and divisorial 
twists. For $\p \in \mathcal{S}(l_0)$, let $\p'$ be a twist of $\p$. 
Let us compare $X_{\p,0}$ and $X_{\p',0}$, where 
$X_{\p,0}$ (resp. $X_{\p',0}$) is the central fiber, that is, the fiber of the map 
$X_{\p} \to \k_0$ (resp. $X_{\p'} \to \k_0$) over $0 \in \k_0$. 
By Proposition 3.1, we 
see that $X_{\p,0}$ and $X_{\p',0}$ are related by a flop  
if and only if the twist is of the 1-st kind or is a small twist of the 2-nd kind. 
For a twist of the 1-st kind, this flop is one of those classified 
in [Na] (Mukai flops of type A, D, $E_6$).  But, for a small twist of the 2-nd kind, 
a new flop appears.  Together with these new flops, we call them Mukai flops 
(cf. Definition 2).     
Let us fix a crepant resolution $\mu_{\p,0}: X_{\p_0,0} \to Y_{\k_0,0}$. 
Take an arbitrary crepant resolution $\mu: X \to Y_{\k_0,0}$. 
Corollary 3.4 states that $X$ is of the form $X_{\p,0}$ for some 
$\p \in \mathcal{S}(l_0)$, and $X$ is related with $X_{\p_0, 0}$ 
by a sequence of Mukai flops.  The proof of Corollary 3.4 goes 
as follows. Let $L' \in \mathrm{Pic}(X)$ be $\mu$-ample  
and let $L \in \mathrm{Pic}(X_{\p_0,0})$ be its proper transform. 
Then $L$ is $\mu_{\p_0,0}$-movable. 
Since $\mathrm{Pic}(X_{\p_0,0})$ is isomorphic to $\mathrm{Pic}(X_{\p_0})$, 
one can lift $L$ to $\mathcal{L} \in \mathrm{Pic}(X_{\p_0})$. 
Then $\mathcal{L}$ is $\mu_{\p_0}$-movable.  By the same argument as in 
Theorem 1.3, we can repeat the flops corresponding to the twists of $\p_0$, 
and finally arrive at $X_{\p_k}$ where the proper transform of $\mathcal{L}$ is $\mu_{\p_k}$-nef. 
Since $L$ is $\mu_{\p_0,0}$-movable, each flops corresponds to a twist of the 1-st kind 
or a small twist of the 2-nd kind.  This implies that the restriction of each flop to the central fibers 
is a Mukai flop.  Then $X_{\p_k,0} = X$ and we get the  desired sequence of Mukai flops 
$X_{\p_0,0} --\to X_{\p_1,0} --\to ... --\to X_{\p_k,0}$. 
    
The author would like to thank 
A. Fujiki for letting him know [Fuj], where he 
learned a key idea for this work. He also thanks 
N. Kawanaka for letting him know 
[Ho] and [Ri].  One of the referees also 
pointed out that there are related works of Howlett and Lehrer 
on $N_W(L_0)/W(L_0)$ (cf. [H-L]).   
Finally he thanks referees for many helpful suggestions.   
\vspace{0.3cm}

{\bf Preliminaries} 

We shall recall some basic facts on nilpotent orbits, 
algebraic groups and birational geometry. Notation 
defined here will be used in the next sections. 
\vspace{0.2cm}

(P.1) {\em Parabolic subgroups and 
Brieskorn-Slodowy diagrams}:  Let $G$ be a complex simple Lie group 
(or a simple algebraic group defined over $\mathbf{C}$) and let $\g$ be its Lie 
algebra. We fix a maximal torus $T$ of $G$ and denote by 
$\h$ its Lie algebra.  Let $\Phi$ be the root system for $\g$ determined by 
$\h$.  The root system 
$\Phi$ has a natural involution $-1$.  With $-1$, one can associate an 
automorphism $\varphi_{\g}$ of $\g$ of order $2$ (cf. [Hu], 14.3). 
This involution will play an important role in Section 1. 
Let 
$$\g = \h \oplus \bigoplus_{\alpha \in \Phi}\g_{\alpha}$$
be the root space decomposition. 
Let us choose a base $\Delta$ of $\Phi$ and denote by 
$\Phi^+$ the set of positive roots. Then 
$\b = \h \oplus 
\bigoplus_{\alpha \in \Phi^+}\g_{\alpha}$ is a Borel subalgebra 
of $\g$ which contains $\h$.  Let $B$ be the corresponding 
Borel subgroup of $G$. 
Take a subset $I$ of $\Delta$.  Let $\Phi_I$ be the root 
subsystem of $\Phi$ generated by $I$ and put 
$\Phi^-_I := \Phi_I \cap \Phi^-$, where $\Phi^{-}$ is the 
set of negative roots. Then 
$$\p_I := \h \oplus\bigoplus_{\alpha \in \Phi^-_I}\g_{\alpha} 
\oplus \bigoplus_{\alpha \in \Phi^+}\g_{\alpha}$$ 
is a parabolic subalagebra containing $\b$. 
This parabolic subalgebra $\p_I$ is called a standard parabolic 
subalgebra with respect to $I$.  Let $P_I$ be the corresponding 
parabolic subgroup of $G$. By definition, $B \subset P_I$. 
Any parabolic subgroup $P$ of $G$ is conjugate to a 
standard parabolic subgroup $P_I$ for some $I$.  Moreover, if 
two subsets $I$, $I'$ of $\Delta$ are different, $P_I$ and 
$P_{I'}$ are not conjugate. Thus, a conjugacy class of parabolic 
subgroups of $G$ is completely determined by $I \subset \Delta$. 
In this paper, to specify the subset $I$ of $\Delta$, we shall 
use the {\em marked} Dynkin diagram. Recall that $\Delta \subset \Phi$ 
defines a Dynkin diagram; each vertex corresponds to a simple root (an 
element of $\Delta$).  Now, if a subset $I$ of $\Delta$ is given, we 
indicate the vertices corresponding to $I$ by white vertices, and 
other vertices by black vertices.  A black vertex is called a {\em marked 
vertex}.  A Dynkin diagram with such a marking is called a {\em 
marked Dynkin diagram}, and a marked Dynkin diagram with only one 
marked vertex is called a {\em single} marked Dynkin diagram. 
Note that the standard parabolic subgroup corresponding to a 
single marked Dynkin diagram is a maximal parabolic subgroup.  
Let $\p$ be a parabolic subalgebra of $\g$ which contains $\h$. 
Let $r(\p)$ (resp. $n(\p)$) be the solvable radical (resp. nilpotent 
radical) of $\p$.  
We put $\k(\p) := r(\p) \cap \h$. Then $$r(\p) = \k(\p) \oplus n(\p).$$  
On the other hand, the Levi factor $l(\p)$ of $\p$ is defined 
as $l(\p) := \g^{\k(\p)}$. Here, 
$\g^{\k(\p)} := \{x \in \g; [x, y] = 0, \forall y \in \k(\p)\}. $ 
Note that $\k(\p)$ is the center of $l(\p)$.  
Then $$\p = l(\p) \oplus n(\p).$$  
If $\p = \p_I$, we have 
$$ \k(\p_I) = \{h \in \h; \alpha(h) = 0, \forall \alpha \in I\}.$$ 
Moreover, we define 
$$ \k(\p_I)^{reg} = \{h \in \k(\p_I); \alpha (h) \ne 0, \forall \alpha \in \Phi\setminus 
\Phi_I\}. $$ Note that $\k(\p_I)^{reg}$ is an open subset of 
$\k(\p_I)$.  
 
Next let us consider the adjoint action of $G$ on $\g$. 
An orbit $\mathcal{O} \subset \g$ for this action is 
called a {\em nilpotent orbit} if $\mathcal{O}$ contains a nilpotent 
element of $\g$. There is a unique nilpotent orbit 
$\mathcal{O}$ such that $n(\p) \cap \mathcal{O}$ is a dense 
open subset of $n(\p)$. Such an orbit is called the {\em Richardson 
orbit} for $\p$. Let $P \subset G$ be the parabolic subgroup 
such that $Lie(P) = \p$. The cotangent bundle $T^*(G/P)$ 
of $G/P$ is isomorphic to the vector bundle $G \times^{P} n(\p)$, 
which is the quotient space of $G \times n(\p)$ by the 
equivalence relation $\sim$. 
Here $(g,x) \sim (g', x')$ if $g' = gp$ and $x' = Ad_{p^{-1}}
(x)$ for some $p \in P$. Let us define a map 
$s: G \times^P n(\p) \to \g$ by $s([g,x]) := Ad_g(x)$. 
Then the image of $s$ coincides with the closure 
$\overline{\mathcal{O}}$ of the 
Richardson orbit for $\p$. The map 
$$ s: T^*(G/P) \to \overline{\mathcal{O}} $$ 
is called the {\em Springer map} for $P$. The Springer map is a 
generically finite, projective morphism.  Thus, the Springer map 
$s$ factorizes as 
$$ T^*(G/P) \stackrel{\pi}\to \tilde{\mathcal O} \to \bar{\mathcal O}.$$ 
The birational map $\pi$ will be a main object in Section 3.  
In particular, 
when $\mathrm{deg}(s) = 1$, we call $s$ the {\em Springer resolution}  
of $\overline{\mathcal{O}}$.  
The following theorem implies that 
$Gr(\p)$ depends only on $\k(\p)$.   
 
\begin{Thm}
$Gr(\p) = 
\overline{G\k(\p)}$. 
\end{Thm}   
   
{\em Proof}. See [B-K], Satz 5.6. 
\vspace{0.12cm}

Every element $x$ of $\g$ can be uniquely written as 
$x = x_n + x_s$ with $x_n$ nilpotent and 
with $x_s$ semi-simple 
such that $[x_n,x_s] = 0$. 
Let $W$ be the Weyl group of $G$ with respect to 
$T$. The set of semi-simple adjoint orbits in $\g$ is 
identified with $\mathfrak{h}/W$. Let 
$\g \to \mathfrak{h}/W$ be the map defined as 
$x \to [\mathcal{O}_{x_s}]$. 
There is 
a direct sum decomposition 
$$ r(\p) = \mathfrak{k}(\p) 
\oplus n(\p), (x \to x_1 + x_2)$$ 
where $n(\p)$ is the 
nil-radical of $\mathfrak{p}$ (cf. [Slo], 
4.3). 
We have a well-defined map 
$$ G \times^P r(\p) \to 
\mathfrak{k}(\p) $$ 
by sending 
$[g, x] \in G \times^P r(\p)$ to 
$x_1 \in \mathfrak{k}(\p)$ and  
there is a commutative diagram 
$$ G \times^P r(\p) \to Gr(\p) $$ 
$$ \downarrow \hspace{1.0cm} \downarrow $$ 
$$ \mathfrak{k}(\p) \to \mathfrak{h}/W, $$
(cf. [Slo], 4.3).  In this paper, we call this a {\bf 
Brieskorn-Slodowy diagram}.  
\vspace{0.2cm}

(P.2) {\em Birational geometry and flops}:             
Let $f: X \to S$ be a projective surjective morphism 
of normal varieties with connected fibers.  
For simplicity, we assume that $S$ is an affine variety 
and $R^1f_*{\mathcal O}_X = 0$.  We are mainly interested in 
the two cases (1) $X$ is a flag variety and $S$ is a point, and 
(2) $S$ is an affine variety with rational singularities and 
$f$ is a resolution of $S$.  We put $N^1(f) := \mathrm{Pic}(X)/f^*
\mathrm{Pic}(S)$ modulo torsion, and $N^1(f )_{\mathbf{R}} = N^1(f) \otimes \mathbf{R}$. 
We sometimes write $N^1(X/S)$ (resp. $N^1(X/S)_{\mathbf{R}}$) 
for $N^1(f)$ (resp. $N^1(f)_{\mathbf{R}}$).  Note that $N^1(f)_{\mathbf{R}}$ 
is a finite dimensional $\mathbf{R}$-vector space. 
In the following sections, we will only consider the case where $S$ is an affine cone (over an 
origin) with a $\mathbf{C}^*$-action with positive weights.  
Then $\mathrm{Pic}(S) = 0$, and $N^1(f)_{\mathbf R} = \mathrm{Pic}(X)\otimes {\mathbf R}$.                
We denote by $N_1(f)$ (or $N_1(X/S)$) the abelian group of numerical  
classes of curves contained in fibers of $f$.  Here two curves $C$ and $C'$ (in some fibers) 
are numerically equivalent if $(L.C) = (L.C')$ for all $L \in \mathrm{Pic}(X)$. 
We define $N_1(f)_{\mathbf{R}} 
:= N_1(f) \otimes_{\mathbf{Z}}\mathbf{R}$. Note that 
$N^1(f)_{\mathbf{R}}$ and $N_1(f)_{\mathbf{R}}$ are dual to each other by 
the intersection form.   
A line bundle $L$ on $X$ is 
called $f$-nef (or nef over $S$) if $(L.C) \geq 0$ for all irreducible proper curves $C$ 
contained in closed fibers of $f$. 
The open convex cone 
$\mathrm{Amp}(f) \subset N^1(f)_{\mathbf{R}}$ (or $\mathrm{Amp}(X/S)$) 
generated by $f$-ample line bundles is called the (relative) {\em ample cone}  
for $f: X \to S$.  Its closure 
$\overline{\mathrm{Amp}}(f) \subset N^1(f)_{\mathbf{R}}$ 
coincides with the closure of the convex cone generated 
by $f$-nef line bundles, and it is called the {\em nef cone} for $f$. 
Note that, when $S$ is affine, $L$ is $f$-ample if and only if $L$ is ample in the 
absolute sense. 
We denote by $\overline{\mathrm{NE}}(f)$ the dual cone of  $\overline{\mathrm{Amp}}(f)$. 
In other words, $$\overline{\mathrm{NE}}(f) = \{z \in N_1(f)_{\mathbf{R}}; 
(L.z) \geq 0, \forall L \in \overline{\mathrm{Amp}}(f)\}. $$     
A line bundle  $L$ on $X$ is called $f$-movable if 
$$ \mathrm{codim} \mathrm{Supp}(\mathrm{coker}(f^*f_*L \to L))  
\geq 2, $$ where $\mathrm{Supp}(F)$ means the support of a 
coherent sheaf $F$.  We denote by 
$\overline{\mathrm{Mov}}(f)$ 
(or $\overline{\mathrm{Mov}}(X/S)$)  the closure of the convex cone 
in $N^1(f)_{\mathbf{R}}$ generated by $f$-movable line bundles.  
Its interior $\mathrm{Mov}(f)$ (or $\mathrm{Mov}(X/S)$) is called 
the (relative) {\em movable cone} for $f$. 
Note that if $f$ is an isomorphism in codimension one, then 
$N^1(f)_{\mathbf R} = \overline{\mathrm{Mov}}(f)$.  
In the next section, we will be concerned with the problem 
of finding all crepant resolutions of an affine variety $S$ with rational 
Gorenstein singularities, when one particular crepant resolution $f: X \to S$ is given.  
Here, a {\em crepant} resolution $f: X \to S$ means a (projective) 
resolution of singularities such that $K_X = f^*K_S$, where 
$K_X$ (resp. $K_S$) is a canonical divisor of $X$ (resp. $K_S$). 
Our strategy for describing other crepant resolutions is as follows. 
Take an arbitrary crepant 
resolution $f': X' \to S$ and let $L'$ be an $f'$-ample line bundle. 
The natural birational map $X --\to X'$ is an isomorphism in 
codimension one because $f$ and $f'$ are both crepant resolutions 
(cf. [K-M], Corollary 3.54).  
Then, one can consider the {\em proper transform} $L \in \mathrm{Pic}(X)$ 
of $L'$.  We want to recover $f': X' \to S$ by using $L$.   
If $L$ is $f$-nef, then $X' = X$; so we assume that $L$ is not $f$-nef. 
Then one can find a birational contraction map $g: X \to \bar{X}$ over $S$ 
such that every curve $C$ contracted by $g$ satisfies $(L.C) < 0$ and 
$g$ is an isomorphism in codimension one.  
We now need a {\em flop} of $g$. A flop of $g$ is a diagram 
$$ X \stackrel{g}\to {\bar X} \stackrel{g'}\leftarrow X_1 $$ such that   
$g'$ is also crepant and the proper transform $L_1 \in \mathrm{Pic}(X_1)$ 
of $L$ is $g'$-nef.  In our situation, we can construct the flop 
explicitly (cf. Section 1); thus, we have another crepant resolution $X_1 \to S$.  
Now, replace $L$ by $L_1$ and continue the same process as above. 
Then we have a sequence of crepant resolutions of $S$: 
$$ X --\to X_1 --\to X_2 --\to ....  $$    
If this sequence terminates at some $X_k$, then $X_k$ is nothing but 
the original $X'$. For details on the birational geometry, see [Ka]. 
\vspace{0.2cm} 

(P.3) {\em The nef cone of a flag variety}:  
Assume that $G$ and $P_I$ are the same as in (P.1).  Then $G/P_I$ is a 
projective manifold. In this case, $\mathrm{Pic}(G/P_I) \cong 
N^1(G/P_I)$.  If $\chi: P_I \to \mathbf{C}^*$ is a character of $P_I$, then 
$P_I$ acts on $\mathbf{C}$ by $\chi$.  Then 
$L_{\chi} := G \times^{P_I} \mathbf{C}$ is a $G$-equivariant line 
bundle on $G/P_I$.  
In this way, we have a natural injective homomorphism 
$$ \phi: \mathrm{Hom}_{\mathrm{alg}.\mathrm{gp}.}(P_I, \mathbf{C}^*) 
\to \mathrm{Pic}(G/P_I), $$ where $\mathrm{Im}(\phi)$ has finite 
index in $\mathrm{Pic}(G/P_I)$. When $G$ is simply connected, this map 
is an isomorphism. Let $L(P_I)$ be the Levi-factor of $P_I$  
corresponding to $l(\p_I)$. Then  
$$\mathrm{Hom}_{\mathrm{alg}.\mathrm{gp}.}(P_I, \mathbf{C}^*) 
\cong \mathrm{Hom}_{\mathrm{alg}.\mathrm{gp}.}(L(P_I), \mathbf{C}^*).$$ 
We define 
$$M(L_{P_I}) : = \mathrm{Hom}_{\mathrm{alg}.\mathrm{gp}.}(L(P_I), \mathbf{C}^*),$$ 
and put $M(L_{P_I})_{\mathbf{R}} := M(L_{P_I}) \otimes \mathbf{R}$.  
Thus, we have 
$$ N^1(G/P_I)_{\mathbf{R}} \cong   
M(L_{P_I})_{\mathbf{R}}.$$  
Recall that $\Delta\setminus I$ corresponds to the set of marked vertices 
$\{v_1, ..., v_{\rho}\}$ of the Dynkin diagram (cf. (P.1)).  Then  
$\rho = \dim N^1(G/P_I)_{\mathbf{R}}$. We set $\Delta_i := \Delta - \{v_i\}$ for 
each $1 \le i \le \rho$.      
The nef cone $\overline{\mathrm{Amp}}(G/P_I)$ is the closed convex cone 
generated by dominant characters $\chi$ of $L(P_I)$ (i.e. $\langle \chi, 
\alpha^{\vee}\rangle   
\geq 0, \forall \alpha \in \Delta$, where $\alpha^{\vee} \in \h$ is the coroot 
corresponding to $\alpha$, and $\chi$ is regarded as an element of 
$\h^*$). Moreover, it is a simplicial cone and each extremal 
ray corresponds to the dominant characters $\chi$ such that 
$\langle\chi, \alpha^{\vee} \rangle = 0, \forall \alpha \in \Delta_i$.  
More geometrically, each extremal ray corresponds to the natural projection  
$p_i: G/P_I \to G/P_{\Delta_i}$ induced by 
the inclusion $P_I \subset P_{\Delta_i}$.        
Now let us consider the Springer map 
$s: T^*(G/P_I) \to \bar{\mathcal{O}}$, where $\mathcal{O}$ 
is the Richardson orbit for $P_I$ (cf. (P.1)). 
Let $$T^*(G/P_I) \stackrel{\pi}\to \tilde{\mathcal{O}} \stackrel{r}\to 
\bar{\mathcal{O}}$$ be the Stein factorization of $s$. 
Then $x_0 := r^{-1}(0)$ consists of one point, and 
$\pi^{-1}(x_0) = G/P_I$. Let $C \subset T^*(G/P_I)$ be a proper curve contained 
in a fiber of $\pi$. Note that $\bar{\mathcal{O}}$ admits a $\mathbf{C}^*$-action with 
positive weights, and this $\mathbf{C}^*$ action extends to an action on  
$T^*(G/P_I)$: $\rho : \mathbf{C}^* \times T^*(G/P_I) \to T^*(G/P_I)$.   Then  
$C_0 := \lim_{t \to 0}\rho(t, C)$ is a curve in $\pi^{-1}(x_0)$ by the 
properness of the relative Hilbert scheme of $\pi$ (cf. [Kol], Claim 1.8.3, Step 5).  
Assume that $L$ is a line bundle on $T^*(G/P_I)$. Then, in order to check if 
$L$ is $\pi$-nef, it is enough to check that $(L.D) \geq 0$ only for 
proper curves $D$ inside $\pi^{-1}(x_0)$. 
Let $p: T^*(G/P_I) \to G/P_I$ be the natural projection map. Then 
we have an isomorphism $p^*: \mathrm{Pic}(G/P_I) \cong  
\mathrm{Pic}(T^*(G/P_I))$. In particular, any line bundle $L$ on 
$T^*(G/P_I)$ has the form $p^*M$ with some $M \in \mathrm{Pic}(G/P_I)$. 
Then, $L$ is $\pi$-nef if and only if $M$ is nef on $G/P_I$.  
As a consequence, there is an identification map 
$$p^*: N^1(G/P_I)_{\mathbf{R}} \cong N^1(\pi)_{\mathbf{R}}$$ such that 
$$p^*(\overline{\mathrm{Amp}}(G/P_I)) = \overline{\mathrm{Amp}}(\pi).$$ 
In particular, the nef cone of $\pi$ is also a simplicial cone. 
\vspace{0.2cm}

\section{} 
Let us consider a  
Brieskorn-Slodowy diagram as in (P.1). 
          
\begin{Lem}
The induced map 
$$ \mu'_{\p}: G \times^P r(\p) \to 
\mathfrak{k}(\p) \times_{\mathfrak{h}/W} 
Gr(\p)$$ 
is a birational projective morphism. In particular, 
$\mathfrak{k}(\p) \times_{\mathfrak{h}/W} 
Gr(\p)$ is irreducible. 
\end{Lem} 

{\em Proof}. The map  
$$G \times^P r(\p) \to G\cdot r(\p)$$ is 
projective. Indeed, it is 
factorized as $G \times^P r(\p) \to G/P \times 
G\cdot r(\p) \to G\cdot r(\p)$, the first 
map is a closed immersion and the second one is a 
projective map because $G/P$ is projective.  
Hence $\mu'_{\p}$ is a projective morphism. 
Let $h \in 
\mathfrak{k}(\p)^{\mathrm{reg}}$ (cf. (P.1)) and 
denote by $\bar{h} \in \mathfrak{h}/W$ its 
image by the map $\mathfrak{k}(P) \to 
\mathfrak{h}/W$. Then the fiber $f^{-1}(\bar{h})$ of 
the map $f: Gr(\p) \to 
\mathfrak{h}/W$ over $\bar{h}$ coincides 
with the semi-simple orbit $G\cdot h$ of $\g$ 
containing $h$. The proof goes as follows.  The centralizer $Z_G(h)$ of $h$ 
is the Levi subgroup $L(P)$ of $P$. Indeed, since $h \in 
\mathfrak{k}(\p)^{\mathrm{reg}}$, $\g^h = l(\p)$. 
This means that $L(P)$ is a subgroup of $Z_G(h)$ 
with finite index. By [Ko], 3.2, Lemma 5, $Z_G(h)$ is connected; 
hence $Z_G(h) = L(P)$. 
Let $U(P)$ be the 
unipotent radical of $P$. Then  
$P\cdot h = U(P)\cdot h$. $U(P)\cdot h$ is closed, and 
its tangent space at $h$ is $h + n(\p)$. 
On the other hand, $h + n(\p)$ is invariant under 
$P$; hence $P\cdot h \subset h + n(\p)$. This implies 
that $P\cdot h = h + n(\p)$. Therefore, $G\cdot (h + n(\p)) 
= G\cdot h$.  By the (Brieskorn-Slodowy) 
diagram, we see that $f^{-1}(\bar{h}) = \cup G\cdot (h+ n(\p)) 
= \cup G\cdot h$, 
where $h$ runs through all elements in the fiber of the map 
$\k(\p) \to \h/W$ over $\bar{h}$. If $h, h' \in \k(\p)$ have 
the same image $\bar{h} \in \h/W$, then $G\cdot h = G\cdot h'$.  
Hence we have $f^{-1}(\bar{h}) = G\cdot h$. 
 
Take a point 
$(h,h') \in \mathfrak{k}(\p)^{\mathrm{reg}} 
\times_{\mathfrak{h}/W} Gr(\p)$. 
Then $h'$ is a semi-simple element 
$G$-conjugate to $h$. 
Fix an element $g_0 \in G$ such that 
$h' = Ad_{g_0}(h)$. 
We have 
$$ (\mu'_{\p})^{-1}(h,h') 
= \{[g,x] \in G \times^P r(\p);\, 
x_1 = h,\, Ad_g(x) = h'\}. $$ 
Since $x = Ad_p(x_1)$ for some $p \in P$ and 
conversely $(Ad_p(x))_1 = x_1$ for any $p \in P$ 
(cf. [Slo], Lemma 2, p.48), we have 
$$ (\mu'_{\p})^{-1}(h,h') = 
\{[g,Ad_p(h)] \in G \times^Pr(\p);\, 
g \in G,\, p \in P,\, Ad_{gp}(h) = h'\} =$$ 
$$ \{[gp,h] \in G \times^Pr(\p);\, 
g \in G,\, p \in P,\, Ad_{gp}(h)  
= h'\} = $$ 
$$ \{[g,h] \in G \times^Pr(\p);\, 
Ad_g(h) = h' \} = $$  
$$ \{[g_0g',h] \in G \times^P r(\p);\, 
g' \in Z_G(h)\} = g_0(Z_G(h)/Z_P(h)). $$  
Since $Z_G(h) = L(P)$, $Z_G(h)/Z_P(h) 
= \{1\}$ and $(\mu'_{\p})^{-1}(h,h')$ 
consists of one point.  Hence, $\mu'_{\p}$ is a birational map 
onto its image. We shall prove that 
$$\mathrm{Im}(\mu'_{\p}) = \mathfrak{k}(\p) \times_{\mathfrak{h}/W} 
Gr(\p).$$ Since $G \times^P r(\p)$ is an affine bundle over  
$G/P$, it is irreducible. Then, the irreducibility of $\mathfrak{k}(\p) \times_{\mathfrak{h}/W} 
Gr(\p)$ follows from this assertion. 
All non-empty fibers of 
$Gr(\p) \to \h/W$ have the same dimension because 
$G \times^P \k(\p) \to \k(\p)$ is an affine bundle. 
So, the map $\mathfrak{k}(\p) \times_{\mathfrak{h}/W} 
Gr(\p) \to \k(\p)$ is equi-dimensional. 
Note that 
$\mathfrak{k}(\p)^{reg} \times_{\mathfrak{h}/W} 
Gr(\p)$ is contained in the closed subset $\mathrm{Im}(\mu'_{\p})$ 
of $\mathfrak{k}(\p) \times_{\mathfrak{h}/W} 
Gr(\p)$.  Hence, $\mathrm{Im}(\mu'_{\p})$ is a unique irreducible 
component of $\mathfrak{k}(\p) \times_{\mathfrak{h}/W} 
Gr(\p)$ of maximal dimension. 
If  $\mathfrak{k}(\p) \times_{\mathfrak{h}/W} 
Gr(\p)$ has an irreducible component different from 
$\mathrm{Im}(\mu'_{\p})$, then its dimension is smaller 
than $\mathrm{Im}(\mu'_{\p})$.  
Let $W'$ be the subgroup of $W$ which stabilizes $\k(\p)$ as 
a set. Then $W'$ acts on $\mathfrak{k}(\p) \times_{\mathfrak{h}/W} 
Gr(\p)$ in such a way that it acts naturally on the first factor and 
it acts trivially on the second factor. 
This $W'$-action must 
preserve $\mathrm{Im}(\mu'_{\p})$.  Assume that $(h,x) \in 
\mathfrak{k}(\p) \times_{\mathfrak{h}/W} Gr(\p)$ is not 
contained in $\mathrm{Im}(\mu'_{\p})$. There is 
an element $h' \in \k(\p)$ such that $h' = w(h)$ with  
some $w \in W'$ and $x \in G(h' + n(\p))$. 
But, $h' \times G(h' + n(\p)) \subset \mathrm{Im}(\mu'_{\p})$ 
and $w^{-1}$ sends $h' \times G(h' + n(\p))$ to 
$h \times G(h' + n(\p))$. As a consequence 
$(h,x) \in w^{-1}(\mathrm{Im}(\mu'_{\p}))$. This is a contradiction. 
Therefore, 
$$\mathrm{Im}(\mu'_{\p}) = \mathfrak{k}(\p) \times_{\mathfrak{h}/W} 
Gr(\p).$$     
Q.E.D. 
\vspace{0.2cm}

In the remainder, we shall write 
$Y_{\k(\p)}$ for the {\em normalization} of 
$\mathfrak{k}(\p) \times_{\mathfrak{h}/W} 
Gr(\p)$, and 
write $X_{\p}$ for 
$G \times^P r(\p)$. 
Then we have a resolution $$\mu_{\p}: X_{\p} \to Y_{\k(\p)}$$  
of $Y_{\k(\p)}$. Let $X_{\p,t}$ (resp. $Y_{\k(\p),t}$) be the 
fiber of the map $X_{\p} \to \k(\p)$ (resp. $Y_{\k(\p)} 
\to \k(\p)$) over $t \in \k(\p)$. Let $\mu_{\p,t}: X_{\p,t} 
\to Y_{\k(\p),t}$ be the map induced by $\mu_{\p}$.   
\vspace{0.2cm} 

\begin{Lem} 
$\mu_{\p}$ is a crepant resolution, which is 
an isomorphism in codimension one. For $t \in 
\k(\p)^{reg}$, $\mu_{\p,t}: X_{\p,t} \to Y_{\k(\p),t}$ is an isomorphism. 
Moreover, for the origin $0 \in \k(\p)$, 
$X_{\p,0} = T^*(G/P)$, $\mu_{\p,0}$ is birational, 
and the Springer map 
$s: T^*(G/P) \to \bar{\mathcal{O}}$ factors 
through $\mu_{\p,0}$ as   
$$ T^*(G/P) \stackrel{\mu_{\p,0}}\to 
Y_{\k(\p),0} \to \bar{\mathcal{O}}.$$ 
\end{Lem}    

{\em Proof}. The second assertion follows from the proof of 
Lemma 1.1. 
By definition, the Springer map $s$ factors through 
$\mu_{\p,0}$. On the other hand, since $\mu_{\p}$ 
is a projective birational morphism and $Y_{\k(\p)}$ 
is a normal variety, $\mu_\p$ has connected fibers by 
Zariski's main theorem; hence $\mu_{\p,0}$ also has 
connected fibers. Since $s$ is a generically finite map, 
we conclude that $\mu_{\p,0}$ is birational.    
These fact imply that, $\mu_{\p,t}$ are birational for 
all $t \in \k(\p)$, and they are isomorphisms for 
general $t$. Therefore,  
$\mu_{\p}$ is an isomorphism in codimension one.
Since $X_{\p,0} = T^*(G/P)$ has trivial canonical 
line bundle, $K_{X_{\p}}$ is $\mu_{\p}$-trivial. 
Then $K_{X_{\p}}$ is the pull-back 
of a line bundle $M$ on $Y_{\k(\p)}$ by $\mu_{\p}$ because 
$Y_{\k(\p)}$ has only rational singularities. But, 
as $Y_{\k(\p)}$ is an affine cone, $M$ is trivial; hence 
$K_{X_{\p}}$ is trivial, and $K_{Y_{\k(\p)}} =( \mu_{\p})_*
K_{X_{\p}}$ is trivial. Q.E.D.   
\vspace{0.2cm}

{\bf Observation 1}. (1) If $\mathrm{deg}(s) > 1$, then 
$\mathfrak{k}(\p) \times_{\mathfrak{h}/W} 
Gr(\p)$ is non-normal. Indeed, the map $\mu'_{\p,0}$ 
(restriction of $\mu'_{\p}$ of Lemma (1.1) to the central 
fibers over $0 \in \k(\p)$), coincides with 
the Springer map $s: T^*(G/P) \to \bar{\mathcal{O}}$, but  
$\mu_{\p,0}$ is birational. 

(2) Let $T^*(G/P) \stackrel{\pi}\to \tilde{\mathcal{O}} 
\to \bar{\mathcal{O}}$ be the Stein factorization of $s$. 
Then $\tilde{\mathcal{O}}$ is the normalization of 
$Y_{\k(\p),0}$ and they coincide at least set-theoretically. 
\vspace{0.2cm} 

We fix a parabolic subgroup $P_0$ of 
$G$ containing $T$, and put $\p_0 := \mathrm{Lie}(P_0)$,  
$\k_0 := \k(\p_0)$ and $l_0 := l(\p)$. We define:  
$$\mathcal{S}(l_0) := \{\p \subset \g; \p: 
\; \mathrm{parabolic} \; s.t. \; 
\h \subset \p, \; \k(\p) = \k_0\}.$$ 
Note that 
$$\mathcal{S}(l_0) = \{\p \subset \g; \p: 
\; \mathrm{parabolic} \; s.t. \; 
\h \subset \p, \; l(\p) = l_0\}.$$ 
The set $\mathcal{S}(l_0)$ is finite. 
Let $\p \in \mathcal{S}(l_0)$ and let $\b$ 
be a Borel subalgebra such that $\h \subset 
\b \subset \p$. Then $\p$ corresponds to a 
marked Dynkin diagram $D$. Take a marked vertex $v$ of 
the Dynkin diagram $D$ and consider the maximal connected 
single marked Dynkin subdiagram $D_v$ of $D$ containing $v$. 
We call $D_v$ the {\em single marked diagram associated with 
$v$}.  When $D_v$ is one of the following, 
we say that $D_v$ (or $v$) {\em is of the 
first kind}, and when $D_v$ does not coincide with any of them, 
we say that $D_v$ (or $v$) {\em is of the second kind}. 
\vspace{0.6cm}

$A_{n-1}$ $(k < n/2)$ 
  
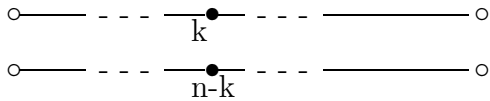
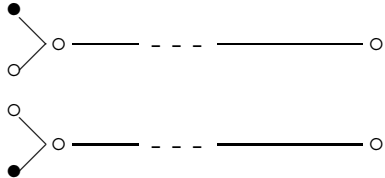
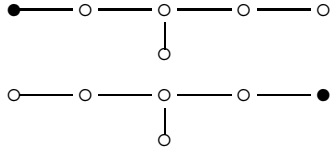
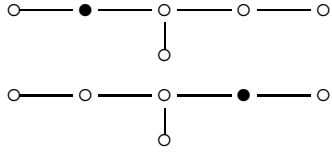
\begin{picture}(300,20)(0,0) 
\put(30,-3){$\circ$}\put(35,0){\line(1,0){25}} 
\put(65,-3.5){- - -}\put(90,0){\line(1,0){15}}
\put(105,-3){$\bullet$}\put(110,0){\line(1,0){10}}
\put(100,-10){k}\put(125,-3.5){- - -}\put(150,0)
{\line(1,0){55}}\put(207,-3){$\circ$}    
\end{picture}  

\begin{picture}(300,20)(0,0) 
\put(30,-3){$\circ$}\put(35,0){\line(1,0){25}} 
\put(65,-3.5){- - -}\put(90,0){\line(1,0){15}}
\put(105,-3){$\bullet$}\put(110,0){\line(1,0){10}}
\put(100,-10){n-k}\put(125,-3.5){- - -}\put(150,0)
{\line(1,0){55}}\put(207,-3){$\circ$}    
\end{picture}
\vspace{0.4cm}

$D_n$ $(n:$ $\mathrm{odd} \geq 5)$  

\begin{picture}(300,20)(0,0) 
\put(30,10){$\bullet$}\put(30,-13){$\circ$}
\put(35,10){\line(1,-1){10}}\put(35,-10){\line(1,1){10}}
\put(47,-3){$\circ$}\put(55,0){\line(1,0){25}}
\put(85,-3.5){- - -}\put(110,0){\line(1,0){55}}\put(167,-3)
{$\circ$}
\end{picture} 
\vspace{0.6cm}

\begin{picture}(300,20)(0,0) 
\put(30,10){$\circ$}\put(30,-13){$\bullet$}
\put(35,10){\line(1,-1){10}}\put(35,-10){\line(1,1){10}}
\put(47,-3){$\circ$}\put(55,0){\line(1,0){25}}
\put(85,-3.5){- - -}\put(110,0){\line(1,0){55}}\put(167,-3)
{$\circ$}
\end{picture}       
\vspace{0.6cm}

$E_{6,I}$:  

\begin{picture}(300,20)
\put(30,-3){$\bullet$}\put(35,0){\line(1,0){20}}
\put(57,-3){$\circ$}\put(65,0){\line(1,0){20}}
\put(87,-3){$\circ$}\put(90,-5){\line(0,-1){10}}
\put(87,-20){$\circ$}\put(95,0){\line(1,0){20}}
\put(117,-3){$\circ$}\put(125,0){\line(1,0){20}}
\put(147,-3){$\circ$} 
\end{picture} 
\vspace{0.4cm}

\begin{picture}(300,20)
\put(30,-3){$\circ$}\put(35,0){\line(1,0){20}}
\put(57,-3){$\circ$}\put(65,0){\line(1,0){20}}
\put(87,-3){$\circ$}\put(90,-5){\line(0,-1){10}}
\put(87,-20){$\circ$}\put(95,0){\line(1,0){20}}
\put(117,-3){$\circ$}\put(125,0){\line(1,0){20}}
\put(147,-3){$\bullet$} 
\end{picture}
\vspace{0.4cm}

$E_{6,II}$: 

\begin{picture}(300,20)
\put(30,-3){$\circ$}\put(35,0){\line(1,0){20}}
\put(57,-3){$\bullet$}\put(65,0){\line(1,0){20}}
\put(87,-3){$\circ$}\put(90,-5){\line(0,-1){10}}
\put(87,-20){$\circ$}\put(95,0){\line(1,0){20}}
\put(117,-3){$\circ$}\put(125,0){\line(1,0){20}}
\put(147,-3){$\circ$} 
\end{picture}
\vspace{0.4cm}

\begin{picture}(300,20)
\put(30,-3){$\circ$}\put(35,0){\line(1,0){20}}
\put(57,-3){$\circ$}\put(65,0){\line(1,0){20}}
\put(87,-3){$\circ$}\put(90,-5){\line(0,-1){10}}
\put(87,-20){$\circ$}\put(95,0){\line(1,0){20}}
\put(117,-3){$\bullet$}\put(125,0){\line(1,0){20}}
\put(147,-3){$\circ$} 
\end{picture}

\vspace{1.0cm}

In the single marked Dynkin diagrams above, two diagrams 
in each type (i.e. $A_{n-1}$, $D_n$, $E_{6,I}$, $E_{6,II}$) 
are called {\em duals}\footnote{ The Weyl group 
$W$ of $\g$ does not contain $-1$ exactly when $\g = A_n (n \geq 2)$, $D_n$ (n: odd) or 
$E_6$ (cf. [Hu], p.71, Exercise 5). This property characterizes the  
Dynkin diagrams in the list.  Moreover, the single marked Dynkin diagrams in the list are  
characterized by the following property. Let $\p_D$ be the parabolic       
subalgebra of $\g$ corresponding to $D$, and let $\varphi_{\g}$ be an automorphism of 
$\g$ determined by $-1$ (cf. (P.1)). 
Then $\varphi_{\g}(\p_D)$ is {\em not}  
conjugate to $\p_D$. If $\varphi_{\g}(\p_D)$ corresponds to a single marked Dynkin 
diagram $D'$, then $D$ and $D'$ are mutually duals. }. 
Let $\bar{D}$ be the marked Dynkin diagram obtained from 
$D$ by making $v$ unmarked. 
Let $\bar{\p}$ be the parabolic subalgebra containing 
$\p$ corresponding to $\bar{D}$. Now let us define a new 
marked Dynkin diagram $D'$ as follows. If $D_v$ is of the 
first kind, we replace $D_v \subset D$ by its dual diagram 
$D_v^*$ to get a new marked Dynkin diagram $D'$. If $D_v$ is 
of the second kind, we define $D' := D$. 
As in (P.1),  the set of unmarked vertices 
of $D$ (resp. $\bar{D}$) defines a subset 
$I \subset \Delta$ (resp. $\bar{I} \subset 
\Delta$). 
By definition, $v \in \bar{I}$. The unmarked vertices of 
$\bar{D}$ define a Dynkin subdiagram, which is decomposed 
into the disjoint sum of the connected component containing $v$ and 
the union of other components. Correspondingly, we have a decomposition 
$\bar{I} = I_v \cup I'_v$ with $v \in I_v$.  
The parabolic subalgebra $\p$ (resp. $\bar{\p}$) coincides 
with the standard parabolic subalgebra $\p_I$ (resp. 
${\p}_{\bar{I}}$). Let $\mathfrak{l}_{\bar{I}}$ be the (standard) 
Levi factor of $\p_{\bar{I}}$. Let $\z(\l_{\bar{I}})$ be the center 
of $\l_{\bar{I}}$. Then $\l_{\bar{I}}/\z(l_{\bar{I}})$ is decomposed 
into the direct sum of simple factors. Now let $\l_{I_v}$ be the 
simple factor corresponding to $I_v$ and let $\l_{I'_v}$ be the 
direct sum of other simple factors. Then 
$$ \mathfrak{l}_{\bar{I}}/\z(\l_{\bar{I}}) = \mathfrak{l}_{I_v} \oplus 
\mathfrak{l}_{I'_v}. $$ 
The marked Dynkin diagram $D_v$ defines a standard parabolic subalgebra 
$\p_v$ of $\mathfrak{l}_{I_v}$. Here let us consider the involution  
$\varphi_{\mathfrak{l}_{I_v}} \in \mathrm{Aut}(\mathfrak{l}_{I_v})$ (cf. (P.1)).       
When $D_v$ is of the first kind, $\varphi_{\mathfrak{l}_{I_v}}(\p_v)$ is 
conjugate to a standard parabolic subalagebra of 
$\mathfrak{l}_{I_v}$ with the dual marked Dynkin diagram 
$D_v^*$ of $D_v$. When $D_v$ is of the second kind, $\varphi_{\mathfrak{l}_{I_v}}(\p_v)$ 
is conjugate to $\p_v$ in $\mathfrak{l}_{I_v}$.  
Let $q: \l_{\bar{I}} \to \l_{\bar{I}}/\z(\l_{\bar{I}})$ be the quotient homomorphism.  
Note that 
$$ \bar{\p} = \mathfrak{l}_{\bar{I}} \oplus n(\bar{\p}), $$ 
$$ \p = q^{-1}(\p_v \oplus \mathfrak{l}_{I'_v}) \oplus n(\bar{\p}). $$ 
Here we define 
$$ \p' = q^{-1}(\varphi_{\mathfrak{l}_{I_v}}(\p_v) \oplus \mathfrak{l}_{I'_v}) \oplus 
n(\bar{\p}). $$ 
Then $\p' \in \mathcal{S}(l_0)$ and $\p'$ is conjugate 
to a standard parabolic subalgebra with the marked Dynkin 
diagram $D'$. {\em This $\p'$ is said to 
be the parabolic subalgebra twisted by $v$}. 

We next define:  
$${\mathcal Res}(Y_{\k_0}) := \{\mathrm{the \: isomorphic \: classes \: 
\: of \: crepant \: projective \: resolutions \: of} \: Y_{\k_0}\}.$$ 
Here we say that two resolutions $\mu : X \to Y_{\k_0}$ and 
$\mu': X' \to Y_{\k_0}$ are isomorphic if there is an isomorphism 
$\phi: X \to X'$ such that $\mu = \mu' \circ \phi$. 
\vspace{0.2cm}

\begin{Thm} 
There is a one-to-one correspondence 
between the two sets $\mathcal{S}(l_0)$ 
and ${\mathcal Res}(Y_{\k_0})$.  Any crepant 
resolution of $Y_{\k_0}$ is an isomorphism in codimension 
one. 
\end{Thm}  

{\em Proof}. By Lemma 1.2, 
$\mu_{\p_0} : X_{\p_0} \to Y_{\k_0}$ is a crepant 
resolution which is an isomorphism in codimension one. Let 
$\mu: X \to Y_{\k_0}$ be an arbitrary crepant resolution. 
Then $X$ and $X_{\p_0}$ are isomorphic in codimension one. 
This implies that $\mu$ is an isomorphism in codimension one. 
This is nothing but the second claim of the theorem. 
Let us consider the first claim.    
Here, the correspondence is given by 
$$ \p \longrightarrow [\mu_{\p}: X_{\p} \to Y_{\k_0}].$$  
In order to check that this correspondence 
is injective, we have to prove that, for $\mu_{\p}: X_{\p} \to Y_{\k_0}$ 
and $\mu_{\p'}: X_{\p'} \to Y_{\k_0}$ with $\p \ne \p'$, 
$\mu_{\p'}^{-1}\circ \mu_{\p}$ is not an isomorphism.  
For $t \in (\k_0)^{reg}$, $X_{\p,t} \cong G \times^P(t + n(\p))$, 
$Y_{\k_0,t} \cong G\cdot t \times \{t\}$, and 
$\mu_{\p,t}$ is an isomorphism defined by 
$$ [g, t+x] \in G \times^P (t+n(\p)) \to (Ad_g(t+x), t) \in 
G\cdot t \times \{t\}.$$  
Since $Z_G(t) = L_0$,  
the map 
$$G \to G \times^P (t+n(\p)) ;\ g \to [g,t]$$ 
induces an isomorphism 
$$ \rho_t: G/L_0 \to G \times^P (t+n(\p)) $$    
In a similar way, we have an isomorphism 
$$ {\rho'}_t: G/L_0 \to G \times^{P'} (t+n(\p')) $$ 
for $t \in (\k_0)^{reg}$.  
The composition of the 
two maps $\rho'_t \circ (\rho_t)^{-1}$ defines 
an isomorphism 
$$ X_{\p,t} \cong X_{{\p'},t} $$ 
for each $t \in (\k_0)^{reg}$ and it 
coincides with $\mu_{\p',t}^{-1}\circ \mu_{\p,t}$. 
Assume that $\mu_{\p'}^{-1}\circ \mu_{\p}$ is a morphism. 
Take $g \in G$ and $p \in P$. We put $g' = gp$. 
Then 
$$\lim_{t \to 0}\rho_t([g]) = [g,0] \in G/P, $$ 
and 
$$\lim_{t \to 0}\rho_t([gp]) = [gp,0] \in G/P. $$ 
Note that $[g,0] = [gp,0]$ in $G/P$. 
On the other hand, 
$$\lim_{t \to 0}\rho'_t([g]) = [g,0] \in G/P', $$ 
and 
$$\lim_{t  \to 0}\rho'_t([gp]) = [gp,0] \in G/P'.$$ 
Since we assume that $\mu_{\p'}^{-1}\circ \mu_{\p}$ 
is a morphism, $[g,0] = [gp,0]$ in $G/P'$. 
But, this implies that $P = P'$, which is a 
contradiction. 
We next prove that the correspondence is surjective. 
Let $\p \in {\mathcal S}(\p_0)$ and let $P$ be the 
corresponding parabolic subgroup of $G$.  
For another $\p' \in \mathcal{S}(l_0)$, $X_{\p'}$ and 
$X_{\p}$ are isomorphic in codimension one; thus, 
$N^1(\mu_{\p'})_{\mathbf{R}}$ and $N^1(\mu_{\p})_{\mathbf{R}}$ 
are naturally identified. 
\vspace{0.2cm}

\begin{Lem} 
(1) $N^1(\mu_{\p})_{\mathbf{R}}$ is identified with 
$N^1(G/P)_{\mathbf{R}}$.  
The nef cone $\overline{\mathrm{Amp}}(\mu_{\p})$ 
is a simplicial polyhedral cone. Each codimension one face 
$F$ of $\overline{\mathrm{Amp}}(\mu_{\p})$ 
corresponds to a marked vertex of the marked Dynkin 
diagram $D$ attached to $P$. 

(2) For each codimension one face $F$, there is an 
element $\p' \in \mathcal{S}(l_0)$ such that 
$$ \overline{\mathrm{Amp}}(\mu_{\p}) \cap 
\overline{\mathrm{Amp}}(\mu_{\p'}) = F.$$ 
\end{Lem} 

{\em Proof}. (1): First of all,  
$\overline{\mathrm{Amp}}(G/P)$ is simplicial by (P.3). Moreover, in (P.3) 
we proved that there is an identification of 
$N^1(\pi)_{\mathbf{R}}$ with 
$N^1(G/P)_{\mathbf{R}}$ which sends  
$\overline{\mathrm{Amp}}(\pi)$ 
to 
$\overline{\mathrm{Amp}}(G/P)$. 
The first statement and the second statement are 
proved in the same manner.  Let us consider the third statement.  
 Let $\bar{D}$ be the marked Dynkin 
diagram which is obtained from $D$ by making a marked 
vertex $v$ unmarked. We then have a parabolic subgroup 
$\bar{P}$ with $P \subset \bar{P}$ corresponding to 
$\bar{D}$. Each codimension one face $F$ of 
$\overline{\mathrm{Amp}}(G/P)$ corresponds to the 
natural surjection $G/P \to G/\bar{P}$. 
We are now going to construct the birational contraction map 
of $X_{\p}$ corresponding to $F$.  
First look at the $\bar{P}$-orbit $\bar{P}r(\p)$ of 
$r(\p)$. Then we can write 
$$\bar{P}r(\p) = r(\bar{\p}) \times L(\bar{P})\cdot r(l(\bar{\p})  
\cap \p),$$ 
where $l(\bar{\p})$ is the Levi factor of $\bar{\p}$, $L(\bar{P})$ is 
the corresponding Levi subgroup, and 
$r(l(\bar{\p}) \cap \p)$ is the solvable radical of 
the parabolic subalgebra $l(\bar{\p}) \cap \p$ of 
$l(\bar{\p})$. 
Then $\k(\p)$ is decomposed as 
$$ \k(\p) = \k(\bar{\p}) \oplus \k(l(\bar{\p}) \cap \p).$$ 
Let $W'$ be the subgroup of the Weyl group $W(l(\bar{\p}))$  
which stabilizes $\k(l(\bar{\p}) \cap \p)$ as a set.  
We then have a map $L(\bar{P})\cdot r(l(\bar{\p}) \cap \p)
\to \k(l(\bar{\p}) \cap \p)/W'$, and hence have 
a map $\bar{P}r(\p) \to \k(l(\bar{\p}) \cap \p)/W'$. 
Now we have a map 
$$ \alpha: \bar{P}\times^P r(\p) \to \bar{P}r(\p) 
\times_{\k(l(\bar{\p}) \cap \p)/W'}\k(l(\bar{\p}) \cap \p).$$ 
Here, the left hand side is isomorphic to 
$$r(\bar{\p}) \times L(\bar{P}) \times^{L(\bar{P}) \cap P} 
r(l(\bar{\p}) \cap \p)$$ and the right hand side is isomorphic 
to $$r(\bar{\p}) \times L(\bar{P})\cdot r(l(\bar{\p}) \cap \p) 
\times_{\k(l(\bar{\p}) \cap \p)/W'}\k(l(\bar{\p}) \cap \p).$$ 
Then, by Lemma (1.1), we see that $\alpha$ is a 
birational map. Since $\alpha$ is a $\bar{P}$-equivariant map, we get  
a birational map:   
$$X_{\p} \to (G \times^{\bar{P}}{\bar{P}}\cdot r(\p)) 
\times_{\k(l(\bar{\p}) \cap \p)/W'}\k(l(\bar{\p}) \cap \p).$$ 
This is the desired birational contraction map.  
 
(2): Let $\p'$ be the parabolic subalgebra twisted by 
$v$. By [Na], Proposition 6.4, we have 
$\bar{P}\cdot n(\p) = \bar{P}\cdot n(\p')$. 
Since $\k(\p) = \k(\p') = \k_0$, we see that 
$$ \bar{P}\cdot r(\p) = \bar{P}\cdot r(\p'). $$ 
Moreover, $$\k(l(\bar{\p}) \cap \p) = \k(l(\bar{\p}) \cap \p').$$ 
Thus, there is a diagram of birational morphisms 
$$ X_{\p} \to 
(G \times^{\bar{P}}{\bar{P}}\cdot r(\p)) 
\times_{\k(l(\bar{\p}) \cap \p)/W'}\k(l(\bar{\p}) \cap \p)
\longleftarrow X_{\p'},$$ and 
$$ \overline{\mathrm{Amp}}(\mu_{\p}) \cap 
\overline{\mathrm{Amp}}(\mu_{\p'}) = F.$$ 

The following corollary implies the 
surjectivity of the correspondence and the proof 
of Theorem 1.3 is completed.  

\begin{Cor} 
Any crepant projective resolution of $Y_{\k_0}$ 
is given by $\mu_{\p}: X_{\p} \to Y_{\k_0}$ 
for some $\p \in \mathcal{S}(l_0)$. 
\end{Cor} 

{\em Proof}. Take an arbitrary crepant projective 
resolution $\mu: X \to Y_{\k_0}$. Fix an $\mu$-ample 
line bundle $L$ on $X$ and denote by $L^{(0)} \in 
\mathrm{Pic}(X_{\p_0})$ its proper transform. 
If $L^{(0)}$ is $\mu_{\p_0}$-nef, then $X$ coincides 
with $X_{\p_0}$. So let us assume that $L^{(0)}$ 
is not $\mu_{\p_0}$-nef. Then there is an extremal 
ray $\mathbf{R}_+[z] \subset \overline{NE}(\mu_{\p_0})$ 
such that $(L^{(0)}.z) < 0$. 
Let $F \subset \overline{\mathrm{Amp}}(\mu_{\p_0})$ 
be the corresponding codimension one face. By the previous 
lemma, one can find $\p_1 \in \mathcal{S}(l_0)$ such 
that  
$$ \overline{\mathrm{Amp}}(\mu_{\p_0}) \cap 
\overline{\mathrm{Amp}}(\mu_{\p_1}) = F.$$ 
We let $L^{(1)} \in \mathrm{Pic}(X_{\p_1})$ be the 
proper transform of $L^{(0)}$ and repeat the same 
procedure. Thus, we get a sequence of flops 
$$ X_{\p_0} --\to X_{\p_1} --\to X_{\p_2} --\to .$$ 
But, since $\mathcal{S}(l_0)$ is a finite set, 
this sequence must terminate by the same argument as 
[Na], Theorem 6.1. 
As a consequence, $X = X_{\p_k}$ for some $k$.  Q.E.D. 
\vspace{0.2cm}

Let $P$ be a parabolic subgroup of $G$ which 
contains $L_0$ as the Levi factor (this is 
equivalent to that $\p \in \mathcal{S}^1(l_0)$).  
By Lemma (1.4), (1), 
$N^1(\mu_{\p})_{\mathbf{R}}$ is identified with 
$N^1(G/P)_{\mathbf{R}}$. By (P.3),  
$N^1(G/P)_{\mathbf{R}}$ is identified with 
$M(L_0)_{\mathbf{R}}$.  Let 
$$ \Phi_P : N^1(\mu_{\p})_{\mathbf{R}} \cong 
M(L_0)_{\mathbf{R}}$$ be the composition of 
these identifications.  
There is a similar identification 
$$ \phi_P : N^1(\pi)_{\mathbf{R}} 
\cong M(L_0)_{\mathbf{R}}.$$ 
The nef cone $\overline{\mathrm{Amp}}(\mu_{\p})$  and 
the (closed) movable cone $\overline{\mathrm{Mov}}(\mu_{\p})$ are regarded 
as cones in $M(L_0)_{\mathbf{R}}$ by $\Phi_P$. 
The nef cone  $\overline{\mathrm{Amp}}(\pi)$ 
and the (closed) movable cone 
$\overline{\mathrm{Mov}}(\pi)$ 
are regarded as cones in $M(L_0)_{\mathbf{R}}$ by $\phi_P$.   
\vspace{0.2cm} 

{\bf Observation 2}. Assume $\p, \p' \in \mathcal{S}(l_0)$. 

(1) $(\Phi_{P'})^{-1}\circ \Phi_P : N^1(\mu_{\p})_{\mathbf{R}} 
\cong N^1(\mu_{\p'})_{\mathbf{R}}$ coincides with the natural 
isomorphism induced by the proper transform by the birational 
map $X_{\p} --\to X_{\p'}$. 

(2) If $X_{\p,0}$ and $X_{\p',0}$ are isomorphic in codimension one, 
then $(\phi_{P'})^{-1}\circ \phi_P : 
N^1(\mu_{\p,0})_{\mathbf{R}} \cong 
N^1(\mu_{\p',0})_{\mathbf{R}}$ coincides with 
the natural 
isomorphism induced by the proper transform by the birational 
map $X_{\p,0} --\to X_{\p',0}$
\vspace{0.2cm} 

{\em Proof}. (1): For $t \in (\k_0)^{reg}$, 
we have an isomorphism $\rho_t : G/L_0 \cong 
X_{\p,t}$ (cf. the first part of the proof of Theorem 1.3). 
Similarly, we have an isomorphism 
$\rho'_t: G/L_0 \cong X_{\p'_t}$.  Let $L_{\chi} \in 
\mathrm{Pic}(X_{\p})$ be the line bundle associated with the 
character $\chi : L_0 \to \mathbf{C}^*$. Then, 
$$(\rho_t)^*L_{\chi}\vert_{X_{\p,t}}  
= G \times^{L_0}_{\chi} \mathbf{C}. $$    
Similarly, let $L'_{\chi} \in \mathrm{Pic}(X_{\p'})$ be the line 
bundle associated with $\chi$. Then 
$$(\rho'_t)^*L'_{\chi}\vert_{X_{\p',t}}  
= G \times^{L_0}_{\chi} \mathbf{C}. $$ 
Thus we see that $L'_{\chi}$ is the proper transform of $L_{\chi}$.    

(2): In this case, the proper transform is compatible with 
the restriction of $X_{\p}$ (resp. $X_{\p'}$) to 
$X_{\p,0}$ (resp. $X_{\p',0}$). Then the result follows from 
(1). 
\vspace{0.2cm}

\begin{Rque} 
The corollary above shows that, for {\em any} 
$L \in \mathrm{Pic}(X_{\p_0})$, there is 
$\p \in \mathcal{S}(l_0)$ such that the 
proper transform $L' \in \mathrm{Pic}(X_{\p})$ 
of $L$ is $\mu_{\p}$-nef. Thus 
$$M(L_0)_{\mathbf{R}} = 
\bigcup_{\p \in \mathcal{S}(l_0)}\overline{\mathrm{Amp}}
(\mu_{\p}). $$ 
\end{Rque} 
 
Let us consider the involution   
$\varphi_{\g} \in \mathrm{Aut}(\g)$ (cf. (P.1)). For $\p \in \mathcal{S}(l_0)$, 
we put $\p^* := \varphi_{\g}(\p)$. Then $\p^* \in \mathcal{S}(l_0)$. 
Note that $(\p^*)^* = \p$. 

\begin{Prop} 
Let $L \in \mathrm{Pic}(X_{\p})$ be a $\mu_{\p}$-ample 
line bundle and let $L' \in \mathrm{Pic}(X_{\p^*})$ be 
its proper transform. Then $-L'$ is $\mu_{\p^*}$-ample. 
In other words, $$-\overline{\mathrm{Amp}}(\mu_{\p}) 
=  \overline{\mathrm{Amp}}(\mu_{\p^*}). $$ 
\end{Prop} 

{\em Proof}. Let $P$ (resp. $P^*$) be the parabolic 
subgroup of $G$ with $\mathrm{Lie}(P) = \p$ (resp. 
$\mathrm{Lie}(P^*) = \p^*$). Note that $P$ and $P^*$ have 
the common Levi factor $L_0$.  
For $t \in (\k_0)^{\mathrm{reg}}$, 
we define an isomorphism 
$$ \rho_t : G/L_0 \to G \times^P (t + n(\p)),$$ 
as in the first part of the proof of Theorem 1.3. 
Denote by the same letter $\varphi$ the automorphism 
of $G$ induced by $\varphi \in Aut(\g)$. 
Note that $\varphi(P) = P^*$ and $P \cap P^* = L_0$. 
The automorphism $\varphi$ 
induces an automorphism $$\bar{\varphi}: G/L_0 \to G/L_0.$$ 
Then 
$$\rho_t \circ \bar{\varphi} \circ (\rho_t)^{-1} : 
G \times^P (t + n(\p)) \to G \times^P (t + n(\p)) $$ 
induces a birational involution 
$$ \sigma : X_{\p} --\to X_{\p}. $$ 
As in the proof of Theorem 1.3, $\sigma$ is not a 
morphism because $P \ne P^*$. 
On the other hand, there is an isomorphism 
$$ \tau: X_{\p} \to X_{\p^*} $$ 
defined by $\tau([g,x]) = [\varphi(g), -{\varphi}(x)]$. 
Then $$\mu_{\p}^{-1}\circ \mu_{\p^*}: X_{\p^*} --\to 
X_{\p}$$ coincides with the composition $\sigma \circ \tau^{-1}$. 
$\mathrm{Pic}(G/L_0) \otimes \mathbf{R}$ is identified 
with $\mathrm{Hom}_{\mathrm{alg}.\mathrm{gp}.}(L_0, \mathbf{C}^*)
\otimes \mathbf{R}$. 
$\mathrm{Hom}_{\mathrm{alg}.\mathrm{gp}.}(L_0, \mathbf{C}^*)$ 
is contained in 
$\mathrm{Hom}_{\mathrm{alg}.\mathrm{gp}.}(T, \mathbf{C}^*)$, 
and $\varphi$ acts on the latter space by $-1$. Therefore,  
$\bar{\varphi}^* \in 
\mathrm{Aut}(\mathrm{Pic}(G/L_0)\otimes \mathbf{R})$ 
is $-1$.           
This means that the proper transform of $L \in 
\mathrm{Pic}(X_{\p}) \otimes \mathbf{R}$ by $\sigma$ is 
$-L$; hence, 
$-L'$ is $\mu_{\p^*}$-ample.  

\begin{Rque} 
(1) Assume $\p \in \mathcal{S}(l_0)$. Assume that the 
Springer map $T^*(G/P) \to \bar{\mathcal{O}}$ is a 
resolution. Then, $\tilde{\mathcal{O}}$ is the normalization of 
$\bar{\mathcal{O}}$. 
Let $\mathcal{S}^1(l_0; \p)$ be the subset of $\mathcal{S}(l_0)$ 
consisting of the parabolic subalgebras $\q$ which are obtained from 
$\p$ by the sequence of twists by marked vertices $v$ of the first kind. 
By Observation 1, $\tilde{\mathcal{O}}$ is the normalization of $Y_{\k_0,0}$. 
By abuse of notation, we denote by the same $\mu_{\q,0}$ the 
Stein factorization $X_{\q,0} \to \tilde{\mathcal{O}}$ of 
$\mu_{\q,0}: X_{\q,0} \to Y_{\k_0,0}$. Then, by [Na], Theorem 6.1, 
$\{\mu_{\q,0}\}_{\q \in 
\mathcal{S}^1(l_0;\p)}$ is nothing but the set of all crepant 
resolutions of $\tilde{\mathcal{O}}$.  
This implies that    
$$ \overline{\mathrm{Mov}}(\pi) 
= \bigcup_{\q \in \mathcal{S}^1(l_0; \p)}\overline{\mathrm{Amp}}(\mu_{\q}).$$  
\vspace{0.2cm} 

(2) When $\p_0$ is a Borel subalgebra $\b_0$ of $\g$, 
$l_0 = \h$. In this case, 
$\mathcal{S}(\h)$ is the set of all Borel subalgebras 
$\b$ such that $\h \subset \b$. The cardinality of 
$\mathcal{S}(\h)$ coincides with the order of the 
Weyl group $W$. In this case, 
$X_{\b} = G \times^{B} \b$, and $Y_{\k_0} = 
\g \times_{\h/W} \h$.  
\end{Rque}   

\begin{Exam} 
Let $\mathcal{O} \subset \g$ be the Richardson orbit for $P_0$. 
Assume that the Springer map $T^*(G/P_0) \to \bar{\mathcal{O}}$ 
is not birational. Let $T^*(G/P_0) \to \tilde{\mathcal{O}}$ 
be the Stein factorization of the Springer map.  
Theorem 1.3 is useful to describe crepant resolutions of 
$\tilde{\mathcal{O}}$. For example, put $\g = so(8,\mathbf{C})$ 
and fix a Cartan subalgebra $\h$ and a Borel subalgebra $\b$ containing $\h$. 
Consider the standard parabolic subalgebra $\p_0$ corresponding 
to the marked Dynkin diagram 

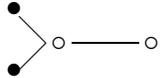
\begin{picture}(120,20)(0,0) 
\put(30,10){$\bullet$}\put(30,-13){$\bullet$}
\put(35,10){\line(1,-1){10}}\put(35,-10){\line(1,1){10}}
\put(47,-3){$\circ$}\put(55,0){\line(1,0){25}}
\put(82,-3){$\circ$}
\end{picture} 
\vspace{0.4cm}

In this case, the Springer map $T^*(G/P_0) \to \bar{\mathcal{O}}$ 
has degree 2. 
In the marked Dynkin diagram above, there are two marked 
vertices. For each marked vertex, one can twist $\p_0$ to 
get a new parabolic subalgebra $\p_1$ or $\p_2$. 
The marked Dynkin diagrams of $\p_1$ and 
$\p_2$ are respectively given by: 

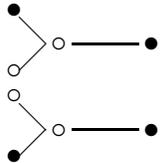
\begin{picture}(120,20)(0,0) 
\put(30,10){$\bullet$}\put(30,-13){$\circ$}
\put(35,10){\line(1,-1){10}}\put(35,-10){\line(1,1){10}}
\put(47,-3){$\circ$}\put(55,0){\line(1,0){25}}
\put(82,-3){$\bullet$}
\end{picture} 
\vspace{0.4cm}

\begin{picture}(120,20)(0,0) 
\put(30,10){$\circ$}\put(30,-13){$\bullet$}
\put(35,10){\line(1,-1){10}}\put(35,-10){\line(1,1){10}}
\put(47,-3){$\circ$}\put(55,0){\line(1,0){25}}
\put(82,-3){$\bullet$}
\end{picture} 
\vspace{0.4cm}

Since $\tilde{\mathcal{O}}$ is the normalization of $Y_{\k_0,0}$,   
the three varieties $X_{\p_0,0}$, $X_{\p_1,0}$ and 
$X_{\p_2,0}$ are crepant resolutions of $\tilde{\mathcal{O}}$. 
But these are not enough; actually, there are 
three more crepant resolutions. In fact, 
$$\mathcal{S}(l_0) = \{\p_0, \p_1, \p_2, (\p_0)^*, 
(\p_1)^*, (\p_2)^*\},$$ and so three other crepant 
resolutions of $\tilde{\mathcal{O}}$ are $X_{(\p_0)^*,0}$, 
$X_{(\p_1)^*,0}$ and $X_{(\p_2)^*,0}$. 
Moreover, $M(L_0)_{\mathbf{R}} = 
\overline{\mathrm{Mov}}
(\pi_0)$ and 
it is divided into the 
union of six ample cones: 

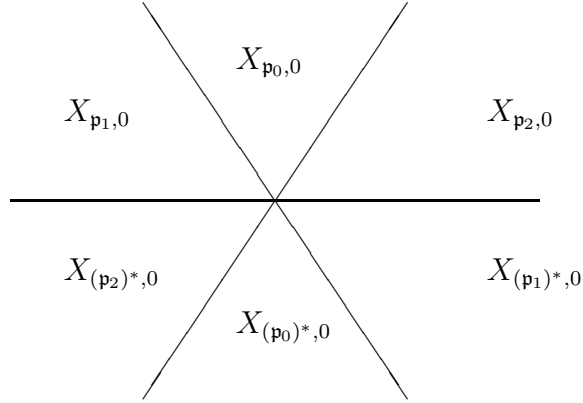
\begin{picture}(300,100)(0,0) 
\put(150,0){\line(2,3){50}}\put(150,0){\line(1,0){100}}
\put(150,0){\line(2,-3){50}}\put(150,0){\line(-2,-3){50}}
\put(150,0){\line(-1,0){100}}\put(150,0){\line(-2,3){50}}
\put(135,50){$X_{\p_0,0}$}\put(230,30){$X_{\p_2,0}$}
\put(230,-30){$X_{(\p_1)^*,0}$}\put(135,-50){$X_{(\p_0)^*,0}$}
\put(70,-30){$X_{(\p_2)^*,0}$}\put(70,30){$X_{\p_1,0}$}
\end{picture} 
\vspace{3.0cm}
\end{Exam}

\section{} 
Let $W$ be the 
Weyl group for $\g$ with respect to $\h$. 
Let  $N_W(L_0) \subset W$ be the normalizer of $L_0$. 
Note that $N_W(L_0) = \{w \in W; w(\k_0) = \k_0\}$. 
For $w \in N_W(L_0)$ and $\chi \in M(L_0)$, we 
define $w\chi \in M(L_0)$ as $w\chi (g) = \chi(w^{-1}gw)$, 
$g \in L_0$. In this way, $N_W(L_0)$ acts on $M(L_0)_{\mathbf{R}}$. 
Note that the Weyl group $W(L_0)$ of $L_0$ is a subgroup of  $N_W(L_0)$ and 
it consists of the elements acting trivially on 
$M(L_0)_{\mathbf{R}}$.   
\vspace{0.2cm}
 
On the other hand, for $\p \in \mathcal{S}(l_0)$, $N_W(L_0)$ 
acts on $N^1(\mu_{\p})_{\mathbf{R}}$ in the following way.    
Since $N_W(L_0)$ acts on $\k_0$, an element 
$w \in N_W(L_0)$ acts 
on $G\cdot r(\p) \times_{\h/W} \k_0$ by 
$id \times w$; hence it acts on the 
normalization $Y_{\k_0}$ of  
$G\cdot r(\p) \times_{\h/W} \k_0$. 
The biregular automorphism of $Y_{\k_0}$ thus defined, induces 
a {\em birational} automorphism $w_{\p}: X_{\p} --\to X_{\p}$. 
Then we have an isomorphism 
$w_{\p}: N^1(\mu_{\p})_{\mathbf{R}} \cong N^1(\mu_{\p})_{\mathbf{R}}$ 
as the push-forward $(w_{p})_*$.  
The following lemma tells us that these two actions are the same 
under the identification $N^1(\mu_{\p})_{\mathbf{R}} \cong 
M(L_0)_{\mathbf{R}}$. 
\vspace{0.2cm}

\begin{Lem} 
(i) $w_{\p}$ coincides with $w$ under the identification 
$N^1(\mu_{\p})_{\mathbf{R}} \cong 
M(L_0)_{\mathbf{R}}$. 

(ii) $w(\overline{\mathrm{Amp}}(\mu_{\p})) = 
\overline{\mathrm{Amp}}(\mu_{w(\p)})$.  
\end{Lem} 

{\em Proof}.           
Note that $W$ acts on the set of parabolic subalgebras of $\g$ 
containing $\h$. Moreover, $N_W(L_0)$ acts on the set 
$\mathcal{S}(l_0)$. In particular, $w(\p) \in 
\mathcal{S}(l_0)$ for $w \in N_W(L_0)$. We define an 
isomorphism $$\phi_w: X_{\p} \to X_{w(\p)}$$ by 
$$ \phi_w ([g,x]) := [gw^{-1}, Ad_w(x)], $$ 
where $[g,x] \in G \times^{P}r(\p)$ and 
$[gw^{-1}, Ad_w(x)] \in G \times^{w(P)} r(w(\p))$. 
The isomorphism $\phi_w$ fits into the commutative diagram  
$$ X_{\p} \stackrel{\phi_w}\to X_{w(\p)}$$ 
$$ \downarrow \hspace{1.0cm} \downarrow $$ 
$$ Y_{\k_0} \stackrel{w}\to Y_{\k_0}. $$  
Restrict this diagram over $t \in (\k_0)^{reg}$. 
Then the vertical maps both induce isomorphisms and 
the horizontal maps induce an isomorphism 
$G/L_0 (= X_{\p,t}) \to G/L_0 (= X_{\p,t})$ given by 
$gL_0 \to gw^{-1}L_0$. Let us choose $\chi \in 
M(L_0)$. Then we get a commutative 
diagram of line bundles over $G/L_0$: 
$$ G \times^{L_0}_{\chi}\mathbf{C} \to 
G \times^{L_0}_{w\chi}\mathbf{C}$$ 
$$ \downarrow \hspace{1.0cm} \downarrow $$ 
$$ G/L_0 \longrightarrow G/L_0. $$  
The left hand side of the horizontal map on the top 
is $L_{\chi}\vert_{X_{\p,t}}$, and the right hand side is 
$L_{w\chi}\vert_{X_{\p,t}}$. This shows that 
$$w_{\p}([L_{\chi}]) = [L_{w\chi}].$$ 
Since $\Phi_P([L_{\chi}]) = \chi$ 
(resp, $\Phi_P([L_{w\chi}]) = w\chi$) by Observation 2, 
we have proved (i). The second claim (ii) is clear from the 
fact that $(\phi_w)_*$ sends $\overline{\mathrm{Amp}}(\mu_{\p})$ 
to $\overline{\mathrm{Amp}}(\mu_{w(\p)})$.   
\vspace{0.2cm}    

\begin{Prop} 
The set $\mathcal{S}(l_0)$ contains exactly $N \cdot \sharp(N_W(L_0)/W(L_0))$ elements, where 
$N$ is the number of the conjugacy classes of parabolic subalgebras 
contained in $\mathcal{S}(l_0)$. 
\end{Prop} 

{\em Proof}.  
 Take two conjugate elements $\p, \p' \in \mathcal{S}(l_0)$.  Then, 
there is an element $w \in W$ 
such that $\p = w(\p')$. Since $r(\p)$ and $r(\p')$ are conjugate 
by $w$, $r(\p) \cap \h$ and $r(\p') \cap \h$ are also conjugate 
by $w$. This means that $w$ sends $\k_0 = r(\p) \cap \h$ to 
$\k_0 = r(\p') \cap \h$, and $w \in N_W(L_0)$.        
We next show that, if  $w(\p) = \p$ for an element $\p \in \mathcal{S}(l_0)$,  then 
$w \in W(L_0)$.   
Let $U$ be the unipotent radical of $P$. 
Then one can write $P = U\cdot L_0$. Now we suppose that 
$w$ is represented by an 
element of the normalizer group $N_G(T)$.  Since $w(P) = P$ 
and $N_G(P) = P$, $w \in P$.    
Let us write $w = u\cdot l$ 
with $u \in U$ and $l \in L_0$. By assumption,  
$w(L_0) = L_0$. This means that $u(L_0) = L_0$.  
Since any two Levi subgroups of $P$ are conjugate by a 
unique element of $U$ (cf. [Bo], 14.19),  we have $u = 1$, 
which implies that $w \in W(L_0)$.   
\vspace{0.2cm}

Let $\mathcal{O} \subset \g$ be the Richardson orbit for 
$P_0$. Assume that the Springer map $s_0: T^*(G/P_0) \to \bar{\mathcal O}$ 
is a birational map. Let $\tilde{\mathcal O}$ be the normalization of $\mathcal{O}$.  
Then $s_0$ factorizes as $T^*(G/P_0) \stackrel{\pi_0}\to \tilde{O} \to \bar{\mathcal O}$.  
The following proposition says that, 
the movable cone $\mathrm{Mov}(\pi_0)$  
is a {\em fundamental domain} for 
the $N_W(L_0)/W(L_0)$-action.  

\begin{Prop} 
Assume that the Springer map $T^*(G/P_0) \to \bar{\mathcal{O}}$ 
is birational.  Define 
$\mathcal{S}^1(l_0)$ to be the subset of $\mathcal{S}(l_0)$ consisting of the 
parabolic subalgebras $\p$ which are obtained from $\p_0$ by the sequence 
of twists of the 1-st kind.  

(i)  The (closed) movable cone $\overline{\mathrm{Mov}}(\pi_0)$ 
coincides with 
$\bigcup_{\q \in \mathcal{S}^1(l_0)}\overline{\mathrm{Amp}}(\mu_{\q})$.

(ii) For any $\p \in \mathcal{S}(l_0)$, there is an element $w \in 
N_W(L_0)$ such that $w(\p) \in \mathcal{S}^1(l_0)$. 

(iii) For any non-zero element $w \in N_W(L_0)/W(L_0)$,  
$$ w(\mathrm{Mov}(\pi_0)) \cap  
\overline{\mathrm{Mov}}(\pi_0) = \emptyset. $$ 
\end{Prop}  

{\em Proof}.  (i) is nothing but Remark 1.8, (1). 

(ii) By the definition of $\mathcal{S}^1(l_0)$, there is an element 
$\p' \in \mathcal{S}^1(l_0)$ such that $\p'$ is conjugate to 
$\p$. Then (ii) follows from the first argument of the proof of 
Proposition 2.2.  

(iii) We shall prove that all elements 
of $\mathcal{S}^1(l_0)$ are not conjugate to each other. 
Suppose that $\p, \p' \in 
\mathcal{S}^1(l_0)$ are conjugate to each other. 
Let us consider the diagram 
$$ X_{\p} \stackrel{\mu_{\p}}\to Y_{\k_0} 
\stackrel{\mu_{\p'}}\leftarrow X_{\p'}. $$ 
Restrict the diagram over $0 \in \k_0$ to get  
$$ X_{\p,0} \stackrel{\mu_{\p,0}}\to Y_{\k_0,0} 
\stackrel{\mu_{\p',0}}\leftarrow X_{\p',0}. $$ 
 Since the Springer map (for $P_0$) has degree 1, $Y_{\k_0,0}$ and 
$\bar{\mathcal{O}}$ are birational and they have the same 
normalization, say $\tilde{\mathcal{O}}$. The Springer maps  
$T^*(G/P) \to \bar{\mathcal{O}}$ and $T^*(G/P') \to \bar{\mathcal{O}}$ 
also have degree 1.  By assumption, $P$ and $P'$ 
are conjugate; hence $T^*(G/P) = T^*(G/P')$ and the two 
Springer maps (resolutions) are the same.  
Let $f_P: T^*(G/P) \to \tilde{\mathcal{O}}$ 
be the Stein factorization of the Springer map for $P$. 
Then $\mu_{\p,0}$ coincides with the composition of $f_P$ 
and the normalization map $\tilde{\mathcal{O}} \to Y_{\k_0,0}$.  
Similarly, $\mu_{\p',0}$ is the composition of $f_{P'}$ and 
the normalization map $\tilde{\mathcal{O}} \to Y_{\k_0,0}$.       
Since $f_P = f_{P'}$, $\mu_{\p,0} = \mu_{\p',0}$.  In particular, the 
birational map $\mu_{\p',0}^{-1} \circ \mu_{\p,0}$ is an isomorphism.  
We shall prove that $\mu_{\p'}^{-1}\circ \mu_{\p}$ is 
an isomorphism. Let $L$ be a $\mu$-ample line 
bundle on $X_{\p}$ and let  $L' \in \mathrm{Pic}(X_{\p'})$  
be the proper transform of $L$ by      
$\mu_{\p'}^{-1}\circ \mu_{\p}$.  
By Corollary 1.5,  
$X_{\p}$ and $X_{\p'}$ are connected by a sequence of 
birational transformations which are isomorphisms in codimension 
one.  Since $\p, \p' \in \mathcal{S}^1(l_0)$, these 
birational tranformations all come from the twist of the first kind. 
This means that, there is a closed subset $F$ of $X_{\p,0}$ with 
codimension $\geq 2$ such that $\mu_{\p'}^{-1}\circ \mu_{\p}$ 
is an isomorphism at each $x \in X_{\p,0}\setminus F$.  
Hence we have   
$$L'\vert_{X_{\p',0}} \cong 
(\mu_{\p',0}^{-1}\circ \mu_{\p,0})_*(L\vert_{X_{\p,0}}).$$ 
But the right hand side is a $\mu_{\p',0}$-ample line bundle. 
Hence $L'\vert_{X_{\p',0}}$ is $\mu_{\p',0}$-ample. 
This shows that $L'$ is $\mu_{\p'}$-ample. Indeed, by the 
$\mathbf{C}^*$-action of $X_{\p'}$, every proper curve $D$ 
in a fiber of $\mu_{\p'}$ is deformed to a curve inside   
 $X_{\p',0}$; hence $(L', D) > 0$ follows from the 
ampleness of $L'\vert_{X_{\p',0}}$. Therefore,    
$\mu_{\p'}^{-1}\circ \mu_{\p}$ is 
an isomorphism.  Then, by Theorem 1.3,  $\p = \p'$. 

\begin{Rque} 
When the Springer map $T^*(G/P_0) \to \bar{\mathcal{O}}$ 
has degree $> 1$, this proposition is not true.  In Example 1.9, the involution 
$-1$ of the root system $\Phi$ of $\g$ can be realized as an element
$w_0$ of the Weyl group $W$ (cf. the footnote in Section 1).    
Clearly, $w_0 \in N_W(L_0)$. But, $w_0$ acts non-trivially 
on $$M(L_0)_{\mathbf{R}} = 
\bigcup_{\p \in \mathcal{S}^1(l_0)}\overline{\mathrm{Amp}}
(\mu_{\p}).$$ 
\end{Rque}  

\section{} 
Recall that,  in Section 1, we have defined operations called {\em twists} 
for the elements $\p \in \mathcal{S}(l_0)$. They are expressed in terms of 
corresponding marked Dynkin diagrams.  
Twists are divided into those of the 1-st kind 
and those of the 2-nd kind. We shall further divide the twists of the 2-nd kind 
into two classes. 
\vspace{0.2cm}

\begin{Def} 
A single marked Dynkin diagram $D$ (of type A,B,C,D,E,F or G) 
of the 2-nd kind is called {\em small} if $D$ is one of the 
following and, otherwise, $D$ is called {\em divisorial}.  
A twist by a small, single marked Dynkin diagram, is  
called of type (2-s), and a twist by a divisorial, single 
marked Dynkin diagram, is called of type (2-d).   
\end{Def} 
\vspace{0.15cm} 
  
$B_n$ (k : even, $k > (2n+1)/3$) 

\begin{picture}(300,20)(0,0) 
\put(30,-3){$\circ$}\put(35,0){\line(1,0){25}} 
\put(65,-3.5){- - -}\put(90,0){\line(1,0){15}}
\put(105,-3){$\bullet$}\put(110,0){\line(1,0){10}}
\put(100,-10){k}\put(125,-3.5){- - -}\put(150,0)
{\line(1,0){15}}\put(170,-3){$\circ$}\put(175,-3){$\Rightarrow$}
\put(187,-3){$\circ$}    
\end{picture}
\vspace{0.4cm}
 
$C_n$ (k : odd, $k \leq 2n/3$) 

\begin{picture}(300,20)(0,0) 
\put(30,-3){$\circ$}\put(35,0){\line(1,0){25}} 
\put(65,-3.5){- - -}\put(90,0){\line(1,0){15}}
\put(105,-3){$\bullet$}\put(110,0){\line(1,0){10}}
\put(100,-10){k}\put(125,-3.5){- - -}\put(150,0)
{\line(1,0){15}}\put(170,-3){$\circ$}\put(175,-3){$\Leftarrow$}
\put(187,-3){$\circ$}    
\end{picture} 
\vspace{0.4cm}
 
$D_n$ (k : odd, $n-2 \geq k > 2n/3$)  

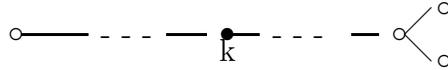
\begin{picture}(300,20)(0,0) 
\put(30,-3){$\circ$}\put(35,0){\line(1,0){25}}
\put(65,-3.5){- - -}\put(90,0){\line(1,0){15}}\put(110,-3)
{$\bullet$}\put(110,-10){k}\put(115,0){\line(1,0){10}}
\put(130,-3.5){- - -} \put(160,0){\line(1,0){10}}
\put(175,-3){$\circ$}
\put(180,0){\line(1,-1){10}}\put(180,0){\line(1,1){10}}
\put(192, 7){$\circ$}\put(192,-13){$\circ$}
\end{picture}
\vspace{0.4cm} 

\begin{Prop} 
Assume that $P \subset G$ is a parabolic subgroup 
conjugate to the standard parabolic subgroup determined 
by a single marked Dynkin diagram $D$ of the 2-nd kind (with 
respect to some Borel subgroup).  Let $\mathcal{O}$ be 
the Richardson orbit for $P$ and let $\pi: T^*(G/P) \to 
\tilde{\mathcal{O}}$ be the Stein factorization of the 
Springer map $T^*(G/P) \to \overline{\mathcal{O}}$. 
Then the birational map $\pi$ is divisorial (i.e. 
$\mathrm{codim}\mathrm{Exc}(\pi) = 1$) if  
$D$ is divisorial, and $\pi$ is small (i.e. 
$\mathrm{codim}\mathrm{Exc}(\pi) \geq 2$) if 
$D$ is small. 
\end{Prop} 

{\em Proof}. 
When $D$ is divisorial, $\pi$ is divisorial by 
[Na], Proposition 5.1. We only have to check that, 
if $D$ is small, then $\mathrm{codim}\mathrm{Exc}(\pi) 
\geq 2$.  First assume that $\g$ is of type $B_n$ (resp.  
$D_n$).  Let $V$ be a $2n + 1$ (resp. $2n$)-dimensional  
$\mathbf{C}$-vector space equipped with a non-degenerate 
symmetric bilinear form $\langle \;, \; \rangle$.  We may assume that 
$G = SO(V)$. A parabolic subgroup $P$ of $G$ is described 
as the stabilizer group of an isotropic flag $F := \{F_i\}_{1 \leq i \leq s}$ 
of $V$.  Here an isotropic flag means a flag such that 
$F_i^{\perp} = F_{s-i}$ for $1 \leq i \leq s$.
The parabolic subgroups $P$ determined by small single marked 
Dynkin diagrams are stabilizer groups of the isotropic flags of 
(flag) type $(k, 2n-2k+1, k)$ (resp. $(k, 2n-2k, k)$). 
When $k$ is even with $k > (2n+1)/3$ 
(resp. $k$ is odd with $k > 2n/3$),  the Richardson orbit 
$\mathcal{O}$ for $P$ has Jordan type 
$[3^{2n-2k+1}, 2^{3k-2n-2}, 1^2]$ (resp. 
$[3^{2n-2k}, 2^{3k-2n-1}, 1^2]$). Moreover, the Springer 
map $s: T^*(G/P) \to \overline{\mathcal{O}}$ has 
degree 2 (see for example, [Na], Section 4). 
When $\overline{\mathcal{O}}$ has no codimension 2 orbits, 
$\mathrm{codim}\mathrm{Exc}(\pi) \geq 2$ by [Na 1], Cor. 1.5.  
Thus, we only have to consider the case where 
$\overline{\mathcal{O}}$ has a codimension 2 orbit, say 
$\mathcal{O}'$. One can check, by using the dimension 
formula ([C-M], Cor. 6.1.4), that this is the case exactly when 
$k = (2n+2)/3$ (resp. $k = (2n+1)/3$). In both cases 
($B_n$ and $D_n$), $\mathcal{O}$ has Jordan type 
$[3^{k-1}, 1^2]$ and $\mathcal{O}'$ has Jordan type $[3^{k-2}, 2^2, 1]$.  
We shall prove that,  for $x \in \mathcal{O}'$, $s^{-1}(x)$ consists 
of one point; here, $s$ is the Springer map. When we regard 
$x$ as an element of  $\mathrm{End}(V)$,  an element of 
$s^{-1}(x)$ corresponds to the isotropic flag $F := \{F_i\}$ of 
type $(k, 2n-2k + 1,k)$ (resp. $(k, 2n-2k,k)$) such that 
$xF_i \subset F_{i-1}$ for all $i$.  Let $\mathbf{d}$ be the 
Young diagram defined by the Jordan type $[3^{k-2}, 2^2, 1]$. 
There is a basis $\{e(i,j)\}$ of $V$ indexed by  
$\mathbf{d}$ with the following properties (cf. [S-S], p.259, 
see also [C-M], 5.1.) 
\vspace{0.12cm}
 
(i) $\{e(i,j)\}$ is a Jordan basis of $x$, that is, $xe(i,j) = e(i-1,j)$ for 
$(i,j) \in \mathbf{d}$. 
 
(ii) $<e(i,j), e(p,q)> \neq 0$ if and only if $p = d_j -i + 1$ and 
$q = \beta(j)$, where $\beta$ is a permutation of 
$\{1, 2, ..., d^1\}$ which satisfies:   
$\beta^2 = id$, $d_{\beta(j)} = d_j$, and 
$\beta(j) \not\equiv j$ (mod 2)
if $d_j \equiv 0$ (mod 2). One can choose an arbitrary 
$\beta$ within these restrictions.   
\vspace{0.12cm}

Introducing such a basis, one can directly check that 
there is only one isotropic flag with the desired properties\footnote{ 
Let $r: \tilde{\mathcal{O}} \to \overline{\mathcal{O}}$ be the 
covering map of degree 2. Then, this observation shows that $r$ 
is ramified along $r^{-1}(\mathcal{O}')$.}. 
    
Next assume that $\g$ is of type $C_n$. 
Let $V$ be a $2n$-dimensional  
$\mathbf{C}$-vector space equipped with a non-degenerate 
skew-symmetric bilinear form $<\;, \;>$.  We may assume that 
$G = Sp(V)$. Similarly to the cases of $B_n$ and $D_n$, a parabolic subgroup 
$P$ of $G$ is described 
as the stabilizer group of an isotropic flag $F := \{F_i\}_{1 \leq i \leq s}$ 
of $V$. 
The parabolic subgroups $P$ determined by small single marked 
Dynkin diagrams are stabilizer groups of the isotropic flags of 
(flag) type $(k, 2n-2k, k)$. 
When  $k$ is odd with $k \geq 2n/3$,  the Richardson orbit 
$\mathcal{O}$ for $P$ has Jordan type 
$[3^{k-1}, 2^2, 1^{2n-3k-1}]$, and the Springer 
map $s$ has degree 2. As in the cases of $B_n$ and $D_n$, we only have 
to discuss the case where  
$\overline{\mathcal{O}}$ has a codimension 2 orbit   
$\mathcal{O}'$. This is the case exactly when 
$k = (2n-1)/3$.  In this case  
$\mathcal{O}$ has Jordan type 
$[3^{k-1}, 2^2]$ and $\mathcal{O}'$ has Jordan type $[3^{k-1}, 2, 1^2]$. 
Now, by the same argument as in the cases of $B_n$ and $D_n$, we 
can prove that $s^{-1}(x)$ consists of one point.       
Q.E.D.   
\vspace{0.2cm}

\begin{Def}  
Let $P \subset G$ be a parabolic subgroup 
conjugate to a standard parabolic subgroup determined 
by a single marked Dynkin diagram $D$. 
Assume that $D$ is of the 2-nd kind and is small. 
Let $\p'$ be the twist of $\p := Lie(P)$. 
Consider the diagram 
$$ X_{\p} \stackrel{\mu_{\p}}\to Y_{\k(\p)} 
\stackrel{\mu_{\p'}}\leftarrow X_{\p'} $$ 
and restrict it to the fibers over $0 \in 
\k(\p)$. Then we have a diagram 
$$ X_{\p,0} \stackrel{\mu_{\p,0}}\to Y_{\k(\p),0} 
\stackrel{\mu_{\p',0}}\leftarrow 
X_{\p',0}.$$ 
Here $X_{\p,0} = T^*(G/P)$ and 
$X_{\p',0} = T^*(G/P')$. 
Let $T^*(G/P) \to \tilde{\mathcal{O}} 
\to \overline{\mathcal{O}}$ be the Stein 
factorization of the Springer map. 
Then $\tilde{\mathcal{O}}$ coincides with 
the normalization of $Y_{\k(\p),0}$. Thus,  
the birational map $\mu_{\p,0}$ (resp. 
$\mu_{\p',0}$) factors through $\tilde{\mathcal{O}}$ 
and we get a diagram 
$$ T^*(G/P) \to \tilde{\mathcal{O}} \leftarrow 
T^*(G/P'). $$ 
By Proposition 3.1, this diagram consists of small 
birational maps. So, this diagram is a flop 
because the original diagram 
$$ X_{\p} \stackrel{\mu_{\p}}\to Y_{\k(\p)} 
\stackrel{\mu_{\p'}}\leftarrow X_{\p'} $$ 
is a flop by Lemma 1.4.   
We call this the Mukai flop of type $B_{n,k}$, 
$C_{n,k}$ or $D_{n,k}$ respectively when the single 
marked Dynkin diagram $D$ is of type   
$B_{n,k}$, $C_{n,k}$ or $D_{n,k}$. 
In [Na], we have defined Mukai flops of type 
$A_{n,k}$, $D_n$ (n:odd), $E_{6,I}$ and 
$E_{6,II}$. Together with these, we call them 
Mukai flops.  
\end{Def} 

\begin{Rque} 
For the Mukai flops of type 
$B_{n,k}$, $C_{n,k}$ and $D_{n,k}$, 
the parabolic subgroups $P$ and $P'$ 
are $G$-conjugate. Hence, there is a natural  
$\overline{\mathcal{O}}$-isomorphism 
$T^*(G/P) \cong T^*(G/P')$. This isomorphism 
induces a {\em non-trivial} 
covering transformation of 
$\tilde{\mathcal{O}} \to \overline{\mathcal{O}}$. 
\end{Rque}

\begin{Exam}
The Mukai flop of type $B_{n,n}$ with $n$ 
even, coincides with the Mukai flop of 
type $D_{n+1}$. In fact, let $G^+_{iso}(n+1, 2n+2)$ and 
$G^{-}_{iso}(n+1, 2n+2)$ be the two connected 
components of the orthogonal Grassmann variety 
parameterising $n+1$-dimensional isotropic 
subspaces of a $2n+2$-dimensional $\mathbf{C}$-
vector space $V$ with a non-degenerate symmetric 
form $<\;, \;>$. Let $V'$ be a $2n+1$-dimensional vector 
subspace of $V$ such that the restriction of $<\;, \;>$ 
to $V'$ is a non-degenerate symmetric form. 
Let $G_{iso}(n, 2n+1)$ be the orthogonal Grassmann 
variety parameterising $n$-dimensional isotropic 
subspaces of $V'$. For an $n+1$-dimensional 
isotropic subspace $W$ of $V$, $W \cap V'$ is an 
isotropic subspace of $V'$ of dimension $n$. This 
correspondence gives isomorphisms
$G^+_{iso}(n+1, 2n+2) \cong G_{iso}(n, 2n+1)$ and 
$G^{-}_{iso}(n+1, 2n+2) \cong G_{iso}(n,2n+1)$. 
These isomorphisms induce the isomorphisms of cotangent 
bundles 
$$T^*(G^+_{iso}(n+1, 2n+2)) \cong T^*(G_{iso}(n,2n+1))$$  
and 
$$T^*(G^{-}_{iso}(n+1, 2n+2)) \cong T^*(G_{iso}(n,2n+1)).$$ 
Let $T^*(G_{iso}(n,2n+1)) \to \overline{\mathcal{O}}$ 
be the Springer map. As noticed in the proof of Proposition 3.1, 
this has degree 2. Then its Stein factorization \\ 
$T^*(G_{iso}(n,2n+1)) \to \tilde{\mathcal{O}}$ coincides 
with the Springer map for \\ $T^*(G^+_{iso}(n+1, 2n+2))$ or 
$T^*(G^{-}_{iso}(n+1, 2n+2))$. Now the Mukai flop of type 
$B_{n,n}$ with $n$ even:  
$$ T^*(G_{iso}(n,2n+1)) \to \tilde{\mathcal{O}} 
\leftarrow T^*(G_{iso}(n,2n+1))$$ 
is identified with the Mukai flop of type $D_{n+1}$:  
$$ T^*(G^+_{iso}(n+1, 2n+2)) \to \tilde{\mathcal{O}} 
\leftarrow T^*(G^{-}_{iso}(n+1, 2n+2)).$$
\end{Exam}

Let $\mathcal{S}^*(l_0)$ be the subset of 
$\mathcal{S}(l_0)$ consisting of the parabolic 
subalgebras $\p$ which are obtained from 
$\p_0$ by a sequence of twists of the 1-st kind 
and twists of type (2-s). 
 
\begin{Cor} 
Let $\pi_0: T^*(G/P_0) \to \tilde{\mathcal{O}}$ 
be the Stein factorization of the Springer map 
$T^*(G/P_0) \to \overline{\mathcal{O}}$. 
Then any two (projective) symplectic resolutions of 
$\tilde{\mathcal{O}}$ are connected by a sequence 
of Mukai flops (Definition 2).   
Moreover, as a cone in $M(L_0)_{\mathbf{R}}$, 
one has 
$$ \overline{\mathrm{Mov}}(\pi_0) 
= \bigcup_{\p \in \mathcal{S}^*(l_0)}
\overline{\mathrm{Amp}}(\mu_{\p}).$$  
\end{Cor} 

{\em Proof}. Take a $\pi_0$-movable line bundle $L$ on 
$T^*(G/P_0)$. We shall prove that 
$$[L] \in \bigcup_{\p \in \mathcal{S}^*(l_0)} 
\overline{\mathrm{Amp}}(\mu_{\p}).$$ 
We may assume that $L$ is not $\pi_0$-nef.  
By the identification 
$$ N^1(\mu_{\p_0})_{\mathbf{R}} 
\cong N^1(\pi_0)_{\mathbf{R}},$$ 
$\overline{\mathrm{Amp}}(\mu_{\p_0})$ and 
$\overline{\mathrm{Amp}}(\pi_0)$ 
are identified. There is an extremal ray 
$\mathbf{R}_+[z] \subset \overline{NE}(\pi_0)$ 
such that $(L.z) < 0$. Let $F \subset 
\overline{\mathrm{Amp}}(\pi_0)$ be the 
corresponding codimension one face. This $F$ can be regarded as 
a codimension one face of $\overline{\mathrm{Amp}}(\mu_{\p_0})$. 
By Lemma 1.4, $F \subset \overline{\mathrm{Amp}}(\mu_{\p_0})$ 
determines a marked vertex $v$ of the marked Dynkin diagram 
$D$ attached to $P_0$. Let $\bar{D}$ be the marked Dynkin diagram 
which is obtained from $D$ by making $v$ unmarked. We then have a 
parabolic subgroup $\bar{P}$ with $P_0 \subset \bar{P}$ corresponding 
to $\bar{D}$. Notation being the same as in the proof of 
Lemma 1.4, the birational morphism determined by $F$ is given by 
$$\phi_F: X_{\p_0} \to (G \times^{\bar{P}}{\bar{P}}\cdot r(\p_0)) 
\times_{\k(l(\bar{\p}) \cap \p_0)/W'}\k(l(\bar{\p}) \cap \p_0).$$   
Note that $\phi_F$ is a morphism over $\k_0$. Restrict 
$\phi_F$ to the fibers over $0 \in \k_0$. Then we have a 
morphism 
$$ \phi_{F,0}: T^*(G/P_0) \to 
G \times^{\bar{P}}{\bar{P}}\cdot n(\p_0). $$ 
The birational morphism $\bar{\phi}_F$ determined by $F \subset 
\overline{\mathrm{Amp}}(\pi_0)$ 
coincides with the Stein factorization of 
$\phi_{F,0}$. 
Since $L$ is $\pi_0$-movable, $\bar{\phi_F}$ must be 
a small birational map, i.e., $\mathrm{codim}(\mathrm{Exc}(\bar{\phi}_F)) 
\geq 2$. On the other hand, 
by Proposition 3.1 and by the 
argument of the proof of [Na], Proposition 6.4, (iii),  
we see that $\bar{\phi}_F$ is a small birational map 
if and only if the single marked Dynkin diagram $D_v$ 
is of the 1-st kind or of the 2-nd kind and small. 
Let $\p_1$ be the twist of $\p_0$ by $v$.  
Then, by Lemma 1.4, we have a flop $X_{\p_0} --\to 
X_{\p_1}$ over $\k_0$. Restrict this flop to the fibers 
over $0 \in \k_0$. Then we have a flop 
$T^*(G/P_0) --\to T^*(G/P_1)$, which is a locally 
trivial family of Mukai flops (cf. [Na], Proposition 
6.4, (iii)). 
Let $L_1$ be the 
proper transform of $L$ by this flop and replace 
$T^*(G/P_0)$ by $T^*(G/P_1)$. The same procedure above 
produces a sequence of flops 
$$ T^*(G/P_0) --\to T^*(G/P_1) --\to T^*(G/P_2) --\to. $$ 
By the same argument as Corollary 1.5, this sequence 
terminates. As a consequence, for some $k$, the proper transform 
$L_k \in \mathrm{Pic}(T^*(G/P_k))$ of $L$ is 
$\pi_k$-nef . 
This implies that 
$$[L] \in \bigcup_{\p \in \mathcal{S}^*(l_0)} 
\overline{\mathrm{Amp}}(\mu_{\p}). $$

{\bf Problem}: (1) {\em Calculate $\sharp \mathcal{S}^*(l_0)$ explicitly.}  

(2) {\em For $\p \in \mathcal{S}(l_0)$, characterize 
an element $\q \in \mathcal{S}^*(l_0)$ such that 
$\mu_{\p,0} = \mu_{\q,0}$.}  

When $\mathrm{deg}(s) = 1$, $\mathcal{S}^*(l_0) = 
\mathcal{S}^1(l_0)$ and any two different elements of the 
set are not conjugate. So, in this case, $\q$ is characterized as a unique element 
of $\mathcal{S}^*(l_0)$ which is conjugate to $\p$. 
However, when $\mathrm{deg}(s) > 1$, $\mathcal{S}^*(l_0)$ may possibly contain  
conjugate, but different elements (cf. Example 1.9).  
The problem is how to characterize $\q$ among the elements of $\mathcal{S}^*(l_0)$ 
which are conjugate to $\p$. 
      
\begin{Rque}  
B. Fu [F, Cor. 5.9] proved that if $s_i: X_i := T^*(G/P_i) \to \bar{\mathcal O}$ ($i = 1,2$) are 
two Springer maps of the same degree, then $X_1$ and $X_2$ 
are connected by a sequence of Mukai flops 
of type $A_{n,k}$, $D_n$ (n:odd), 
$E_{6,I}$ and $E_{6,II}$: 
$Y_j \stackrel{\phi_j}-\to Y_{j+1}$ ($j =0, ..., k-1$) 
with $Y_0 = X_1$ and $Y_k = X_2$.  
But, his result is in a different context from ours. In fact, let  
$\pi_1: X_1 \to  \tilde{\mathcal O}$ be the Stein factorization of $s_1$. 
Then $\phi_j$ are birational maps over $\bar{\mathcal{O}}$, but not 
necessarily over $\tilde{\mathcal{O}}$.  More exactly, each $\phi_j$ may 
possibly induce a non-trivial covering automorphism of 
$\tilde{\mathcal O} \to \bar{\mathcal O}$.  
For example, in Example 1.9, $X_{\p_i,0}$ and $X_{(\p_i)^*,0}$ are  
isomorphic as the varieties over $\bar{\mathcal O}$, but not isomorphic 
as those over $\tilde{\mathcal O}$.   
\end{Rque}

\quad \\
\quad\\

Yoshinori Namikawa \\
Department of Mathematics, 
Graduate School of Science, Osaka University, JAPAN \\
namikawa@math.sci.osaka-u.ac.jp


\begin{thebibliography}{}

\bibitem[B-K]{B-K} Borho, W., Kraft, H.:  
\"{U}ber Bahnen und deren Deformationen bei linearen 
Aktionen reduktiver Gruppen, Comment Math. Helv. 
{\bf 54} (1979), 61-104 

\bibitem[Bo]{Bo} Borel, A.: Linear algebraic groups, second 
enlarged edition, GTM {\bf 126}, Springer, 1991  

\bibitem[C-M]{C-M} Collingwood, D., McGovern, W.: 
Nilpotent orbits in semi-simple Lie algebras, 
van Nostrand Reinhold, Math. Series, 1993 

\bibitem[F]{F} Fu, B.: 
Extremal contractions, stratified Mukai flops and 
Springer maps, Adv. Math. {\bf 213} (2007), 165-182 

\bibitem[Fuj]{Fuj} Fujiki, A.: 
A generalization of an example of Nagata, 
in the proceeding of Algebraic Geometry Symposium, 
Jan. 2000,  Kyushu University, 111-118 

\bibitem[He]{He} Hesselink, W.: 
Polarizations in the 
classical groups, Math. Z. {\bf 160} (1978) 
217-234 

\bibitem[H-L]{H-L} Howlett, R., Lehrer, G.: 
Induced cuspidal representations and generalized 
Hecke rings, Invent. Math. {\bf 58} (1980) 37-64
  
\bibitem[Ho]{Ho} Howlett, R.: Normalizers of parabolic 
subgroups of reflection groups, J. London Math. Soc. 
{\bf 21}  (1980), 62-80  

\bibitem[Hu]{Hu} Humphreys, J.: 
Introduction to Lie algebras and representation 
theory, Graduate Texts in Math. {\bf 9}, Springer 
(1972) 

\bibitem[Ka]{Ka} Kawamata, Y.: 
Crepant blowing-up of 3-dimensional canonical 
singularities and its application to 
degenerations of surfaces, Ann. Math. {\bf 127}, 
93-163 (1988)  

\bibitem[Ko]{Ko} Kostant, B.: 
Lie group representations 
on polynomial rings, Amer. J. Math. {\bf 85} 
(1963) 327-404 

\bibitem[Kol]{Kol} Koll{\'a}r, J.: Rational curves on algebraic 
varieties, Springer, Berlin Heidelberg New York, 1996 

\bibitem[K-M]{K-M} Koll{\'a}r, J., Mori, S.: 
Birational geometry of algebraic varieties, Cambridge University Press,  
1998 

\bibitem[Na]{Na} Namikawa, Y.: Birational geometry 
of symplectic resolutions of nilpotent orbits, 
to appear in Advances Studies in Pure Mathematics {\bf 45}, 
(2006), Moduli Spaces and Arithmetic Geometry (Kyoto, 2004), 
pp. 75-116, see also math.AG/0404072,  math.AG/0408274

\bibitem[Na 1]{Na 1} Namikawa, Y.: 
Deformation theory of singular 
symplectic n-folds, Math. Ann. {\bf 319} (2001) 
597-623

\bibitem[Ri]{Ri} Richardson, R.W.: 
Conjugacy classes of involutions in Coxeter groups.  Bull. Austral. Math. Soc.  
{\bf 26}  (1982), no. 1, 1-15

\bibitem[Slo]{Slo} Slodowy, P.: Simple singularities 
and simple algebraic groups. Lect. Note Math. 
{\bf 815}, Springer-Verlag, 1980

\bibitem[S-S]{S-S} Springer,T.A., Steinberg, R.: 
Conjugacy classes, In: Borel, et. al.: 
seminar on algebraic groups and related finite 
groups, Lect. Note Math. {\bf 131} 167-266, 
Springer Verlag 1970    


\end{thebibliography}
\end{document}